\documentclass[10pt]{amsart} 

\usepackage{amsmath}
\usepackage{amssymb}

\usepackage{graphicx}
\usepackage{xcolor}
\usepackage{multirow}
\usepackage{enumitem}

\newtheorem{theorem}{Theorem}[section]
\theoremstyle{definition}
\newtheorem{example}[theorem]{Example}

\newcommand{\bfu}{\mathbf u}
\newcommand{\bfw}{\mathbf w}
\newcommand{\bfD}{\mathbf D}
\newcommand{\bfR}{\mathbf R}
\newcommand{\dx}{\mathrm d}
\newcommand{\FD}{\text{FD}}
\newcommand{\LSZ}{\text{\tiny LSZ}}
\newcommand{\MWENO}{\text{\tiny MWENO}}
\newcommand{\CWENO}{\text{\tiny CWENO}}

\newcommand{\C}{\text{\tiny C}}

\newcommand{\cfl}{\text{CFL}}

\newcommand{\exact}{\text{exact}}

\title[WENO Approximations to Sharp Propagating Fronts for RD Systems]{A Numerical Study of WENO Approximations to Sharp Propagating Fronts for Reaction-Diffusion Systems}

\author{Jiaxi Gu}
\address{Department of Mathematics $\&$ POSTECH MINDS (Mathematical Institute for Data Science), Pohang University of Science and Technology, Pohang 37673, Korea}
\email{jiaxigu@postech.ac.kr}

\author[2]{Daniel Olmos-Liceaga}
\address{Departmento de Matem\'{a}ticas, Universidad de Sonora, Rosales y Blvd. Encinas, S/N, Hermosillo, Sonora 83000, Mexico}
\email{daniel.olmos@unison.mx}

\author{Jae-Hun Jung}
\address{Department of Mathematics $\&$ POSTECH MINDS (Mathematical Institute for Data Science), Pohang University of Science and Technology, Pohang 37673, Korea}
\email{jung153@postech.ac.kr}

\makeatletter
\@namedef{subjclassname@2020}{\textup{}2020 Mathematics Subject Classification}
\makeatother

\begin{document}

\maketitle

\begin{abstract}
Many reaction-diffusion systems in various applications exhibit traveling wave solutions that evolve on multiple spatio-temporal scales. 
These traveling wave solutions are crucial for understanding the underlying dynamics of the system. 
However, because these systems are generally nonlinear, exact solutions are not available for many cases, making the task of obtaining fast and accurate solutions challenging. 
In this work, we present sixth-order weighted essentially non-oscillatory (WENO) methods within the finite difference framework to solve reaction-diffusion systems. 
High-order finite difference methods, commonly used in the reaction-diffusion research community, tend to become inaccurate and unstable when the solution is highly localized and nonsmooth unless a large number of spatial grid points and small time steps are used. 
The WENO method allows us to use fewer grid points and larger time steps compared to classical finite difference methods. 
Our focus is on solving the reaction-diffusion system for the traveling wave solution with the sharp front, for which the WENO method is particularly well suited.  
Although the WENO method is popular for hyperbolic conservation laws, especially for problems with  discontinuity, it can be adapted for the equations of parabolic type, such as reaction-diffusion systems, to effectively handle sharp wave fronts. 
Thus, we employed the WENO methods specifically developed for equations of parabolic type. We considered various reaction-diffusion equations, including the Fisher's, Zeldovich, Newell-Whitehead-Segel, bistable equations and the Lotka-Volterra competition-diffusion system, all of which yield traveling wave solutions with sharp wave fronts. 
Numerical examples in this work demonstrate that the central WENO method is highly more accurate and efficient than the commonly used finite difference method. 
We also provide an analysis related to the numerical speed of the sharp propagating front in the Newell-Whitehead-Segel equation. 
We found that there exists an optimal value of the CFL number for the given number of spatial grid points, where the numerical speed matches the exact one for the Newell-Whitehead-Segel equation. The overall results confirm that the central WENO method is highly efficient and is recommended for solving reaction-diffusion equations with sharp wave fronts. 
\end{abstract}

\section{Introduction} \label{sec:intro}
Traveling waves form a fundamental class of solutions to reaction-diffusion systems, which correspond to many physical, chemical, and biological phenomena \cite{Epstein,MurrayI,MurrayII}.
In this work, we are interested in the numerical solutions of the traveling waves for the one-dimensional reaction-diffusion equations of the following form
\begin{equation} \label{eq:1d_rd}
 \frac{\partial \bfu}{\partial t} = \bfD \frac{\partial^2 \bfu}{\partial x^2} + \bfR (\bfu),
\end{equation}
where $\bfu(x,t) = (u_1(x,t), \cdots, u_n(x,t))^T$ is a vector of $n$ concentration variables, $\bfD$ denotes a diagonal matrix of the diffusion coefficients and $\bfR(\bfu) = (R_1(\bfu), \cdots, R_n(\bfu))^T$ describes the nonlinear source term including the reaction coefficient $\rho$ (see Section \ref{sec:1d_rd} for details).
The existence of traveling wave solutions depends on the reaction term $\bfR(\bfu)$ and the analytic solutions are not available due to the nonlinearity of the reaction term in general. 
Thus, those solutions are obtained numerically. 
In order to approximate \eqref{eq:1d_rd} numerically, various numerical methods such as finite difference \cite{He,WangZhang}, finite volume \cite{Coudiere,Penland,Zhou}, finite element \cite{Franzone,Lakkis,Tuncer} and Chebyshev pseudospectral \cite{Jung2023,Olmos2015} methods, have been developed.

When these numerical methods are applied to solve reaction-diffusion equations, simulating sharp wave fronts presents a significant challenge and is of great interest. 
For example, the large value of the reaction coefficient $\rho$ in $R(u)$ gives rise to the sharp front of the traveling wave, separating one equilibrium state from the other. 
Due to the presence of such sharp wave fronts, the classical finite difference methods require a large number of spatial grid points and very small time steps for better accuracy. 
If the number of spatial grid points is insufficient and the time step is too large, instability may quickly arise.
One of the most important objectives in calculating reaction-diffusion systems is to accurately capture the speed of the wave fronts. 

Since the sharp front propagates with time, its behavior looks similar to the behavior of shock solutions to  hyperbolic conservation laws. 
Therefore, it is reasonable to apply the weighted essentially non-oscillatory (WENO) method to such reaction-diffusion systems, enabling the sharp front to be well resolved.
The original WENO methods were specifically developed for hyperbolic conservation laws, with all nonlinear weights carefully designed to avoid the spurious oscillations around discontinuities.
Since reaction-diffusion equations are of parabolic type, it is reasonable to adapt the WENO method to handle the parabolic term effectively. 
In \cite{Liu}, Liu et al. constructed the sixth-order finite difference WENO (WENO-LSZ) scheme for the degenerate parabolic equation, which approximates the second derivative term directly by a conservative flux difference.
However, unlike the positive linear weights of WENO schemes for hyperbolic conservation laws \cite{JiangShu,Shu}, the negative linear weights exist so that some special cares, such as the technique in \cite{Shi}, were applied to guarantee the non-oscillatory performance in regions of sharp interfaces.
Following the definition of the smoothness indicators in \cite{JiangShu,Shu} and invoking the mapped function, the resulting nonlinear weights meet the requirement of sixth-order accuracy.
In \cite{Hajipour}, Hajipour and Malek proposed the modified WENO (MWENO) scheme with Z-type nonlinear weights \cite{Borges} and nonstandard Runge–Kutta (NRK) schemes.
Further, the hybrid scheme, based on the spatial MWENO and the temporal NRK schemes, was employed to solve the degenerate parabolic equation numerically.
In \cite{Rathan}, Rathan and Gu proposed the central WENO (CWENO) scheme to avoid the negative linear weights and the splitting technique for computational efficiency. 

In this paper, we apply those WENO methods developed for the equations of parabolic type to various reaction-diffusion systems, including the Fisher's, Zeldovich, Newell–Whitehead–Segel and bistable equations, and the Lotka-Volterra competition-diffusion system, and compare them with the classical centered finite difference method. 
We first show the numerical results of different methods as the number of spatial grid points increases. 
Then we compare the numerical results for large time steps and discuss the speed issue in the Newell–Whitehead–Segel equation. 
The basic property that characterizes the traveling wave solutions is the finite speed of propagation, which is one of the numerical challenges in simulating such reaction-diffusion systems.
One solution is presented for matching the finite speed of propagation well. 

In this paper, we demonstrate that the CWENO method is a promising alternative to classical numerical methods, such as finite difference methods, for accurately and efficiently solving reaction-diffusion equations while matching the propagating speed of the sharp wave front. 
Particularly, through the current work, we find that there could be the optimal value of CFL number for the given number of grid points, with which the numerical and exact speeds agree in the Newell–Whitehead–Segel equation. 
That is, decreasing arbitrarily the CFL number does not necessarily provide the monotonic convergence behavior of the numerical speed to the exact one.  
But instead, there exists an optimal combination of the number of grid points and the CFL number, providing the best approximation to the numerical speed. 
With this optimality, the CFL number does not need to be small for the numerical speed to be close to the exact one. 
The phenomenon in this paper is not found for linear problems and most equations that we considered. 
In particular, this phenomenon occurs for the Newell-Whitehead-Segel equation. 

The rest of this paper is organized as follows. 
Section \ref{sec:1d_rd} gives a brief review of the reaction-diffusion systems to be solved numerically.
In Section \ref{sec:weno}, the sixth-order WENO methods, that are developed particularly to handle the diffusion terms, are presented. 
In Section \ref{sec:nr} we provide the numerical results about the accuracy of the WENO methods and address the time-stepping and speed issues in solving the reaction-diffusion systems with the sharp wave fronts.
The concluding section provides some remark and outlines the future developments.

\section{One-dimensional reaction-diffusion system} \label{sec:1d_rd}
A special type of the solution to \eqref{eq:1d_rd} is the traveling wave solution of the form $\bfu(x,t) = \bfw(x-ct)$ for a particular speed $c$, where $\bfw$ are monotonic, and for $z = x-ct$,
$$ \lim_{z \to -\infty} \bfw(z) = \bfu_-, \quad \lim_{z \to \infty} \bfw(z) = \bfu_+, $$
which means that the solution switches from the equilibrium state $\bfu=\bfu_-$ to the other equilibrium state $\bfu=\bfu_+$.

The scalar reaction-diffusion equation in one dimension has the form
\begin{equation} \label{eq:1d_scalar}
 \frac{\partial u}{\partial t} = D \frac{\partial^2 u}{\partial x^2} + R(u),
\end{equation}
with $D$ the constant diffusion coefficient and $R(u)$ a local reaction kinetics.
If $R(u) = \rho u(1-u)$, \eqref{eq:1d_scalar} is known as the Fisher's equation \cite{Fisher}, where its exact solution is given by the following
\begin{equation} \label{eq:fisher_exact}
 u(x,t) = \frac{1}{\left\{ 1 + \exp \left[ \sqrt{\frac{\rho}{6 D}} \left( x - 5 \sqrt{\frac{\rho D}{6}} t \right) \right] \right\}^2}.
\end{equation}
As shown in the exact solution of the traveling wave, the wave speed is determined by the diffusion coefficient $D$ and the reaction coefficient $\rho$. 
The Zeldovich equation \cite{Gilding} has the source term of $R(u) = \rho u^2(1-u)$.
Its analytic solution is given by the following 
\begin{equation} \label{eq:zeldovich_exact}
 u(x,t) = \frac{1}{1 + \exp \left[ \sqrt{\frac{\rho}{2D}} \left( x - \sqrt{\frac{\rho D}{2}} t \right) \right]}.
\end{equation}
With the reaction term $R(u) = \rho u (1-u^\alpha)$, \eqref{eq:1d_scalar} corresponds to the Newell–Whitehead–Segel equation, whose exact solution \cite{Wang} is given as the following form
\begin{equation} \label{eq:nws_exact}
 u(x,t) = \left\{ \frac{1}{2} \tanh \left[ - \frac{\alpha}{2 \sqrt{2\alpha+4}} \sqrt{\frac{\rho}{D}} \left( x - \frac{\alpha+4}{\sqrt{2\alpha+4}} \sqrt{\rho D} t \right) \right]+ \frac{1}{2} \right\}^{2/\alpha}.
\end{equation}
As shown in the above exact form, the wave speed is determined by the diffusion coefficient $D$, reaction coefficient $\rho$ and the exponent of $\alpha$. 
The bistable equation \cite{Gilding}, with the source term $R(u) = \rho u(1-u)(u-\beta)$ for $0 < \beta < 1$, has an analytic solution of the traveling wave form
\begin{equation} \label{eq:bistable_exact}
 u(x,t) = \frac{1+\beta}{2} + \frac{1-\beta}{2} \tanh \left[ \frac{1-\beta}{4} \sqrt{\frac{2 \rho}{D}} \left( x + (1+\beta) \sqrt{\frac{\rho D}{2}} t \right) \right].
\end{equation}
Again as shown in the exact solution above, the wave speed is determined by the diffusion coefficient $D$, reaction coefficient $\rho$ and the left limit equilibrium state $\beta$.

\begin{figure}[htbp]
\centering
\includegraphics[width=0.45\textwidth]{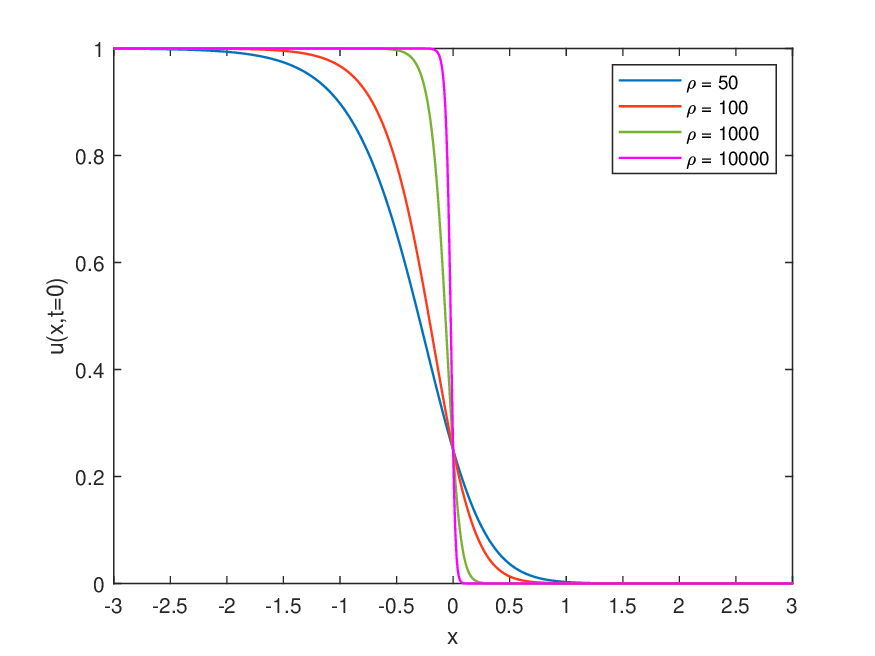}
\includegraphics[width=0.45\textwidth]{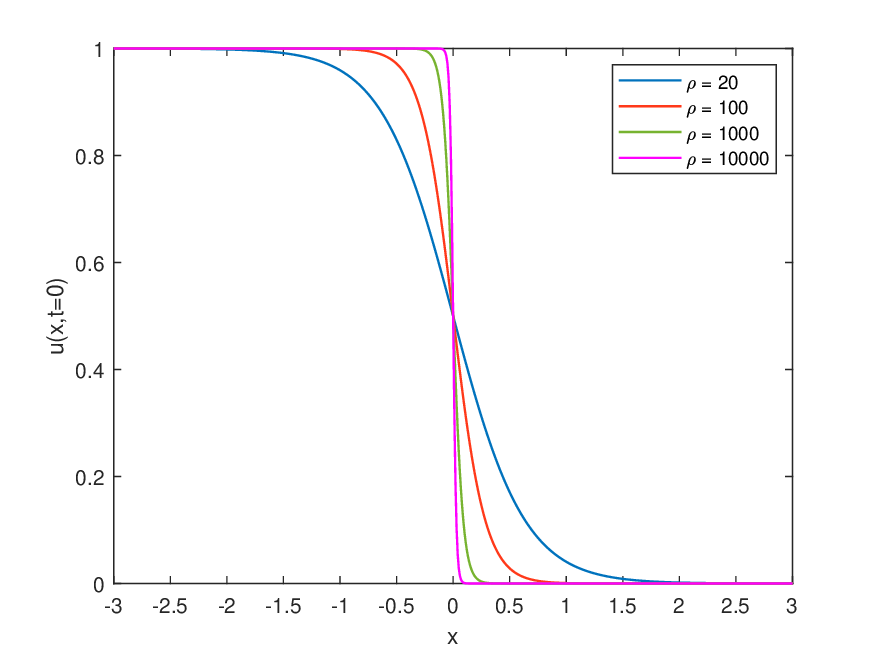}
\includegraphics[width=0.45\textwidth]{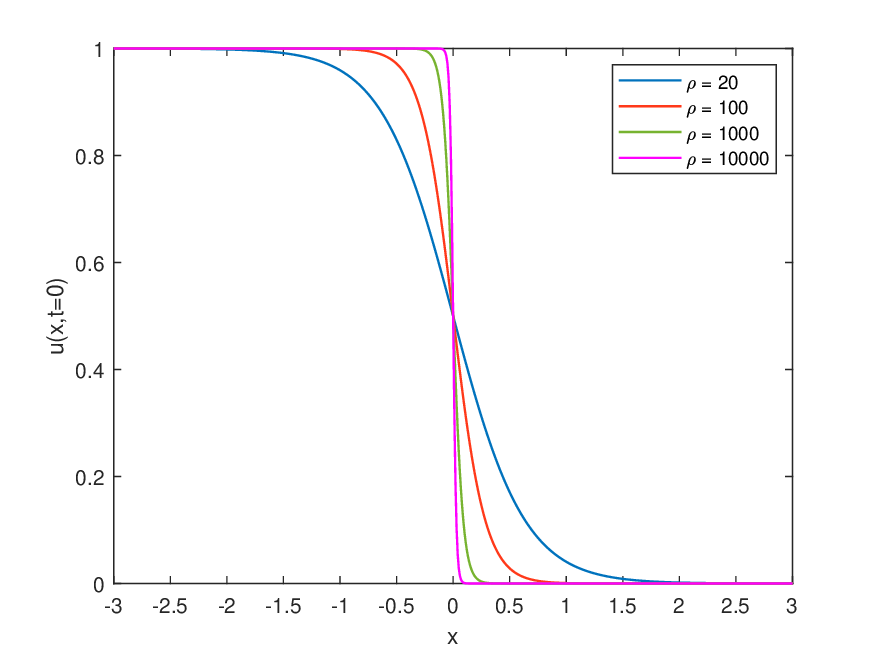}
\includegraphics[width=0.45\textwidth]{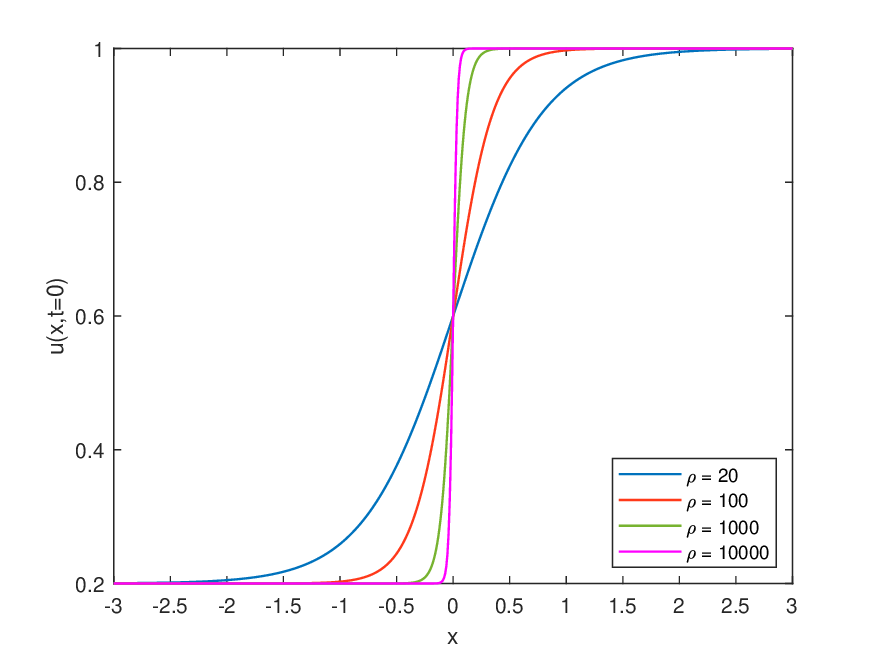}
\caption{Exact solutions to the Fisher's equation (top left), Zeldovich equation (top right), Newell–Whitehead–Segel equation with $\alpha=2$ (bottom left) and bistable equation with $\beta=0.2$ (bottom right) at $t=0$ for different values of $\rho$.}
\label{fig:exact_rho}
\end{figure}

\begin{figure}[htbp]
\centering
\includegraphics[width=0.5\textwidth]{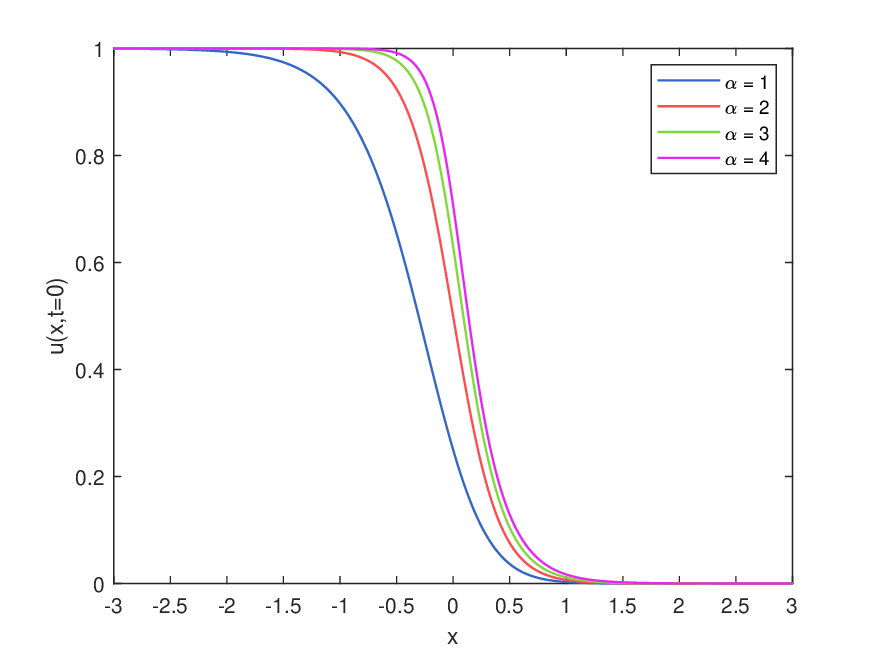}
\vspace{-0.2in}
\caption{Exact solutions \eqref{eq:nws_exact} to the Newell–Whitehead–Segel equation ($\rho=50$) at $t=0$ for different values of $\alpha$.}
\label{fig:nws_alpha}
\end{figure}

\begin{figure}[htbp]
\centering
\includegraphics[width=0.5\textwidth]{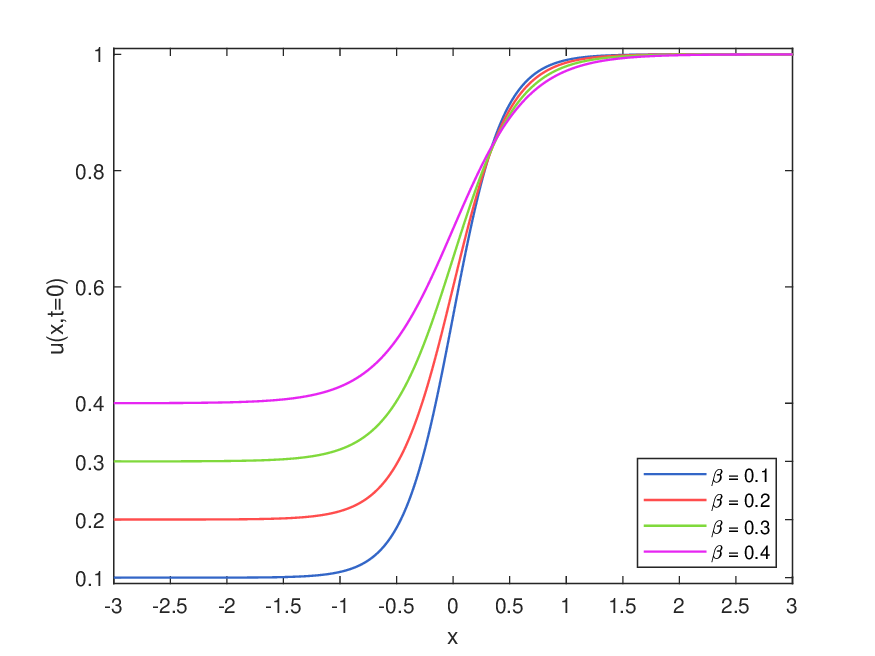}
\vspace{-0.2in}
\caption{Exact solutions \eqref{eq:bistable_exact} to the bistable equation ($\rho=50$) at $t=0$ for different values of $\beta$.}
\label{fig:bistable_beta}
\end{figure}

Figure \ref{fig:exact_rho} shows how the exact solutions to these four equations at $t=0$ behave for various values of $\rho$.
Note that as $\rho$ becomes larger, the solution tends toward a sharper function with the rapid transition around $x=0$ as shown in the figure.
Figure \ref{fig:nws_alpha} presents the behavior of the exact solution \eqref{eq:nws_exact} to the Newell–Whitehead–Segel equation at $t=0$ for various values of $\alpha$, where the solution becomes sharper as $\alpha$ increases.
The appearance of the exact solution \eqref{eq:bistable_exact} to the bistable equation at $t=0$ for various values of $\beta$, is presented in Figure \ref{fig:bistable_beta}.
Note that we have 
$$ \lim_{x \to -\infty} u(x,0) = \beta, $$
at the left boundary.

\begin{figure}[htbp]
\centering
\includegraphics[width=0.45\textwidth]{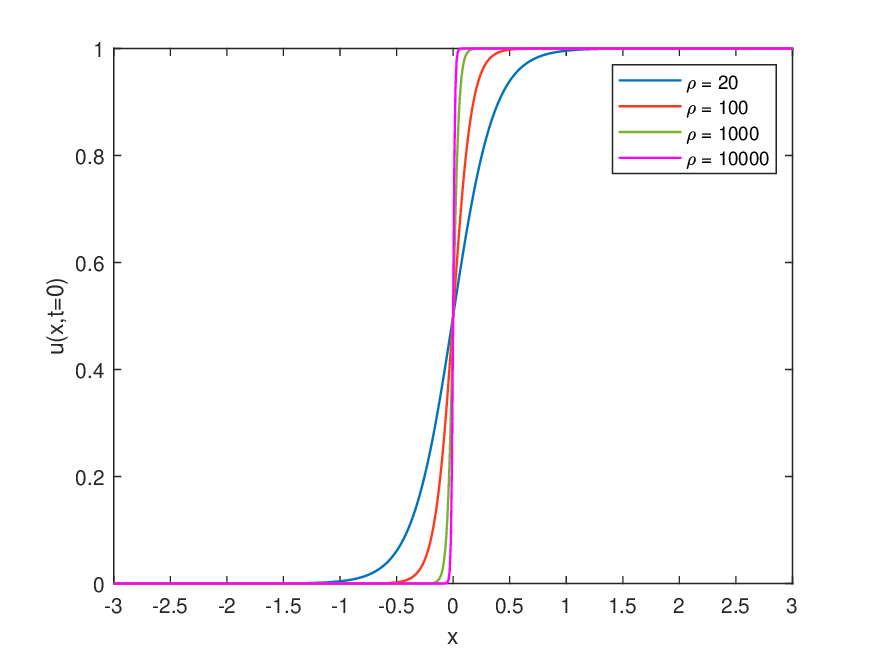}
\includegraphics[width=0.45\textwidth]{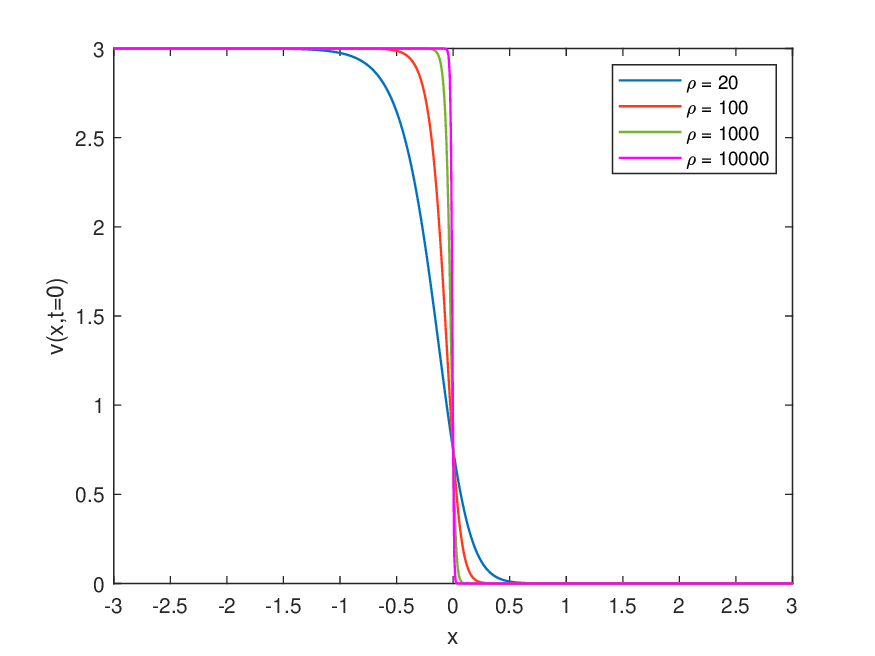}
\caption{Exact solutions $u$ (left) and $v$ (right) to Lotka-Volterra competition-diffusion system \eqref{eq:lotka_volterra} at $t=0$ for different values of $\rho$.}
\label{fig:lotka_volterra_exact_rho}
\end{figure}

The Lotka-Volterra competition-diffusion system is of the form
\begin{equation} \label{eq:lotka_volterra}
\begin{aligned}
 \frac{\partial u}{\partial t} &= D \frac{\partial^2 u}{\partial x^2} + \rho u \left( 1-u-v \right), \\
 \frac{\partial v}{\partial t} &= \frac{D}{3} \frac{\partial^2 v}{\partial x^2} + \rho v \left( 3-4u-v \right),
\end{aligned}    
\end{equation}
whose exact solutions are given as below 
\begin{equation} \label{eq:lotka_volterra_exact}
\begin{aligned}
 u(x,t) &= \frac{1}{2} \left\{ 1 + \tanh \left[ \frac{1}{2} \sqrt{\frac{3\rho}{2 D}} \left( x - \sqrt{\frac{\rho D}{6}} t \right) \right] \right\}, \\
 v(x,t) &= \frac{3}{4} \left\{ 1 - \tanh \left[ \frac{1}{2} \sqrt{\frac{3\rho}{2 D}} \left( x - \sqrt{\frac{\rho D}{6}} t \right) \right] \right\}^2.
\end{aligned}
\end{equation}
Figure \ref{fig:lotka_volterra_exact_rho} shows how the exact solution $(u,v)$ to this reaction-diffusion system at $t=0$ appears for various values of $\rho$.

\section{WENO approximation to reaction-diffusion equation} \label{sec:weno}
In this section, we explain sixth-order finite difference WENO schemes for solving the scalar reaction-diffusion equation \eqref{eq:1d_scalar}.
Given the computational domain $[a, \, b]$, we  apply a uniform grid with $N+1$ points,
$$ x_i = a + i \Delta x, \quad i = 0, \cdots, N, $$
with $\Delta x = (b-a)/N$,  the grid spacing.
Each grid point $x_i$ is also called the cell center for the $i$th cell $I_i = [x_{i-1/2}, \, x_{i+1/2}]$, where the cell boundaries are defined as $x_{i \pm 1/2} = x_i \pm \Delta x/2$.
The semi-discrete form of \eqref{eq:1d_scalar} with respect to $t$, yields
\begin{equation} \label{eq:rd_discrete}
 \frac{du(x_i,t)}{dt} = D \left. \frac{\partial^2 u}{\partial x^2} \right|_{x=x_i} + R(u(x_i, t)).
\end{equation} 

\begin{figure}[htbp]
\centering
\includegraphics[width=\textwidth]{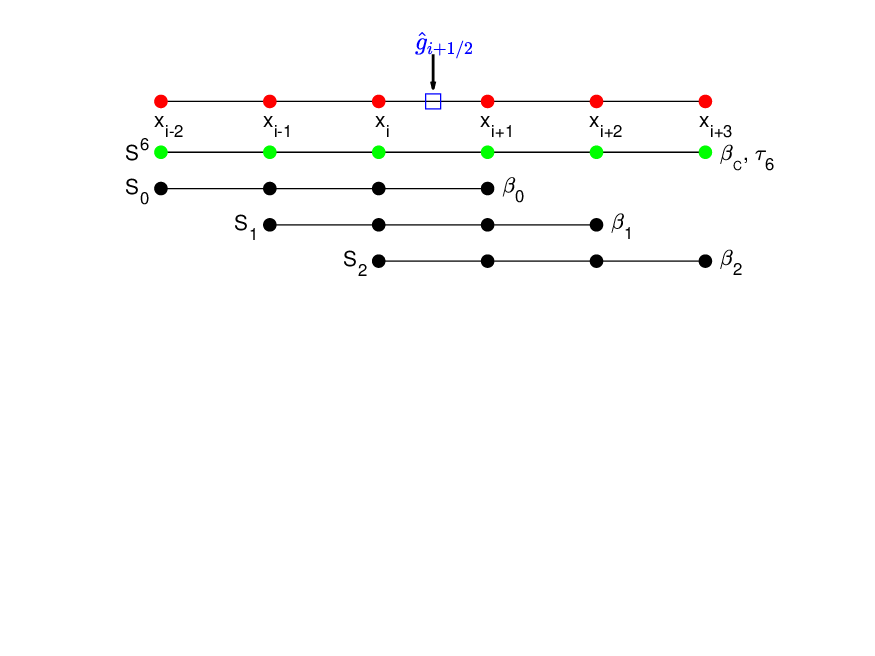}
\vspace{-2.3in}
\caption{The numerical flux $\hat{g}_{i+1/2}$ is constructed on the stencil $S^6 = \{ x_{i-2}, \cdots, x_{i+3} \}$ with six uniform points, as well as three $4$-point substencils $S_0, S_1, S_2$.}
\label{fig:stencil}
\end{figure}

Define the function $h(x)$ implicitly by
$$
   u(x) = \frac{1}{\Delta x^2} \int^{x+\Delta x/2}_{x-\Delta x/2} \left( \, \int^{\eta+\Delta x/2}_{\eta-\Delta x/2} h(\xi) \dx \xi \right) \dx \eta.
$$
Differentiating both sides twice with respect to $x$ gives
$$ 
   \frac{\partial^2 u}{\partial x^2} = \frac{h(x+\Delta x) - 2 h(x) + h(x-\Delta x)}{\Delta x^2}. 
$$
Setting $g(x) = h(x+\Delta x/2) - h(x-\Delta x/2)$, we obtain
\begin{equation} \label{eq:partial_u}
 \left. \frac{\partial^2 u}{\partial x^2} \right|_{x=x_i} = \frac{g_{i+1/2} - g_{i-1/2}}{\Delta x^2}, 
\end{equation} 
where $g_{i \pm 1/2} = g(x_{i \pm 1/2})$.
In order to approximate $g_{i+1/2}$, a polynomial approximation $q(x)$ \cite{Rathan} to $h(x)$ of degree at most $5$ is constructed on the $6$-point stencil $S^6$, as shown in Figure \ref{fig:stencil}. 
The polynomial approximation $p(x)$ to $g(x)$ is obtained by taking the difference of $q(x+\Delta x/2)$ and $q(x-\Delta x/2)$, that is, 
$$
   g(x) = h(x+\Delta x/2) - h(x-\Delta x/2) \approx q(x+\Delta x/2) - q(x-\Delta x/2) = p(x).
$$
Evaluating $p(x)$ at $x = x_{i+1/2}$ yields the numerical flux
\begin{equation} \label{eq:numerical_flux_plus}
 \hat{g}_{i+1/2}^{\FD} = p(x_{i+1/2}) 
                       = - \frac{1}{90} u_{i-2} + \frac{5}{36} u_{i-1} - \frac{49}{36} u_i + \frac{49}{36} u_{i+1} - \frac{5}{36} u_{i+2} + \frac{1}{90} u_{i+3}.   
\end{equation}
The numerical flux $\hat{g}_{i-1/2}^{\FD}$ is obtained by shifting one grid to the left, which gives
\begin{equation} \label{eq:numerical_flux_minus}
 \hat{g}_{i-1/2}^{\FD} = - \frac{1}{90} u_{i-3} + \frac{5}{36} u_{i-2} - \frac{49}{36} u_{i-1} + \frac{49}{36} u_i - \frac{5}{36} u_{i+1} + \frac{1}{90} u_{i+2}.
\end{equation}
Applying Taylor expansions to $\hat{g}_{i \pm 1/2}^{\FD}$ \eqref{eq:numerical_flux_plus} and \eqref{eq:numerical_flux_minus}, we have
\begin{align}
 \hat{g}_{i+1/2}^{\FD} &= g_{i+1/2} + \frac{1}{560} h_i^{(7)} \Delta x^7 + O(\Delta x^8), \label{eq:numerical_flux_plus_Taylor} \\
 \hat{g}_{i-1/2}^{\FD} &= g_{i-1/2} + \frac{1}{560} h_i^{(7)} \Delta x^7 + O(\Delta x^8). \label{eq:numerical_flux_minus_Taylor}
\end{align}
Replacing $g_{i \pm 1/2}$ in \eqref{eq:partial_u} with \eqref{eq:numerical_flux_plus_Taylor} and \eqref{eq:numerical_flux_minus_Taylor} yields the sixth-order centered finite difference approximation
\begin{equation} \label{eq:FD_scheme}
 \left. \frac{\partial^2 u}{\partial x^2} \right|_{x=x_i} = \frac{\hat{g}_{i+1/2}^{\FD} - \hat{g}_{i-1/2}^{\FD}}{\Delta x^2} + O(\Delta x^6).
\end{equation}

A similar approach can be used to obtain a polynomial $p_k(x)$ of degree at most $2$ on each $4$-point substentil $S_k = \{ x_{i-2+k}, \cdots, x_{i+1+k} \}$ with $k = 0,1,2$, resulting in the numerical flux $\hat{g}^k_{i+1/2}$,
\begin{equation} \label{eq:numerical_flux_substencil}
\begin{aligned}
 \hat{g}^0_{i+1/2} &= \frac{1}{12} u_{i-2} - \frac{1}{4} u_{i-1} - \frac{3}{4} u_i + \frac{11}{12} u_{i+1}, \\
 \hat{g}^1_{i+1/2} &= \frac{1}{12} u_{i-1} - \frac{5}{4} u_i + \frac{5}{4} u_{i+1} - \frac{1}{12} u_{i+2}, \\
 \hat{g}^2_{i+1/2} &= -\frac{11}{12} u_i + \frac{3}{4} u_{i+1} + \frac{1}{4} u_{i+2} - \frac{1}{12} u_{i+3}.
\end{aligned}
\end{equation}
In addition, there exist linear weights $d_0 = d_2 = - \frac{2}{15}$ and $d_1 = \frac{19}{15}$ such that
$$ \hat{g}_{i+1/2}^{\FD} = \sum_{k=0}^2 d_k \hat{g}^k_{i+1/2}, $$
and an index shift by $-1$ gives the corresponding relation between $\hat{g}_{i-1/2}^{\FD}$ and $\hat{g}^k_{i-1/2}$.
The WENO procedure cannot be applied directly to the numerical fluxes $\hat{g}^k_{i+1/2}$ \eqref{eq:numerical_flux_substencil} as the existence of negative linear weights, which may lead to the blow-up numerical solutions \cite{Shi}.
So we apply the splitting technique in \cite{Shi} to the linear negative weights,
\begin{gather*}
 \tilde{\gamma}^+_k = \frac{1}{2} \left( d_k + 3 | d_k | \right), \quad \tilde{\gamma}^-_k = \tilde{\gamma}^+_k - d_k, \\
 \sigma^{\pm} = \sum_{k=0}^2 \tilde{\gamma}_k^{\pm}, \quad \gamma_k^{\pm} = \tilde{\gamma}_k^{\pm} / \sigma^{\pm},
\end{gather*}
for $k=0,1,2$.
Then the positive/negative sums $\sigma^{\pm}$, and the linear positive/negative weights $\gamma_k^{\pm}$ are given by
\begin{equation} \label{eq:gamma}
\begin{aligned}
 \sigma^+ &= \frac{42}{15}, & \gamma^+_0 &= \frac{1}{21}, & \gamma^+_1 &= \frac{19}{21}, & \gamma^+_2 &= \frac{1}{21}; \\
 \sigma^- &= \frac{27}{15}, & \gamma^-_0 &= \frac{4}{27}, & \gamma^-_1 &= \frac{19}{27}, & \gamma^+_2 &= \frac{4}{27},
\end{aligned}
\end{equation}
which satisfy 
\begin{equation} \label{eq:linear_weight_relation}
 d_k = \sigma^+ \gamma_k^+ - \sigma^- \gamma_k^-.
\end{equation}
In \cite{Liu}, the nonlinear positive and negative weights $\omega^{\pm}_k$ are defined as
\begin{equation} \label{eq:weights_Liu_pm}
 \omega^{\pm}_k = \frac{\alpha^{\pm}_k}{\sum^2_{l=0} \alpha^{\pm}_l}, \quad \alpha^{\pm}_k = \frac{\gamma^{\pm}_k}{(\beta_k + \epsilon)^2},\quad 
  k=0,1,2,
\end{equation}
where the smoothness indicators are
\begin{equation} \label{eq:smooth_indicator_Liu}
\begin{aligned}
 \beta_0 &= \frac{13}{12} \left( u_{i-2} - 3 u_{i-1} + 3 u_i - u_{i+1} \right)^2 + \frac{1}{4} \left( u_{i-2} - 5 u_{i-1} + 7 u_i - 3 u_{i+1} \right)^2, \\
 \beta_1 &= \frac{13}{12} \left( u_{i-1} - 3 u_i + 3 u_{i+1} - u_{i+2} \right)^2 + \frac{1}{4} \left( u_{i-1} - u_i - u_{i+1} + u_{i+2} \right)^2,\\
 \beta_2 &= \frac{13}{12} \left( u_i - 3 u_{i+1} + 3 u_{i+2} - u_{i+3} \right)^2 + \frac{1}{4} \left( -3 u_i + 7 u_{i+1} - 5 u_{i+2} + u_{i+3} \right)^2,
\end{aligned}
\end{equation}
and $\epsilon>0$, e.g. $\epsilon = 10^{-6}$, avoids the denominator becoming zero. 
Based on the relation \eqref{eq:linear_weight_relation} for the linear weights, the nonlinear weights are defined by
\begin{equation} \label{eq:weights}
 \omega^*_k = \sigma^+ \omega^+_k - \sigma^- \omega^-_k.
\end{equation}
Then the mapped functions are employed to increase the accuracy, 
$$
 g_k(\omega) = \frac{\omega (d_k + d_k^2 - 3 d_k \omega + \omega^2)}{d_k^2 + \omega (1 - 2d_k)}, \quad k = 0,1,2.
$$
The final nonlinear weights are defined as
$$
 \omega^{\LSZ}_k = \frac{\alpha_k}{\sum^2_{l=0} \alpha_l}, \quad
 \alpha_k = g_k(\omega^*_k), \quad k=0,1,2.
$$
In \cite{Hajipour}, the MWENO scheme was proposed by using the Z-type nonlinear weights \cite{Borges} with the global smoothness indicator $\tau = |\beta_0 - \beta_2|$.
The nonlinear positive and negative weights are defined as 
\begin{equation} \label{eq:weights_MWENO_pm}
 \omega^{\pm}_k = \frac{\alpha^{\pm}_k}{\sum^2_{l=0} \alpha^{\pm}_l}, \quad \alpha^{\pm}_k = \gamma^{\pm}_k \left( 1 + \left( \frac{\tau}{\beta_k + \epsilon} \right)^2 \right), \quad k=0,1,2,
\end{equation}
with $\gamma^{\pm}_k$ in \eqref{eq:gamma} and $\epsilon = 10^{-30}$.
As defined in \eqref{eq:weights}, the MWENO nonlinear weights are 
$$
   \omega^{\MWENO}_k = \sigma^+ \omega^+_k - \sigma^- \omega^-_k.
$$
The WENO numerical flux for the above two WENO techniques is
$$
   \hat{g}_{i+1/2} = \sum_{k=0}^2 \omega_k \hat{g}^k_{i+1/2},
$$
where $\hat{g}^k_{i+1/2}, \, k=0,1,2$ are given by \eqref{eq:numerical_flux_substencil}.

In order to avoid the negative linear weights, the central numerical flux $\hat{g}^\C_{i+1/2}$
\begin{equation} \label{eq:numerical_flux_central}
 \hat{g}^\C_{i+1/2} = - \frac{3}{40} u_{i-2} + \frac{11}{24} u_{i-1} - 2u_i + 2u_{i+1} - \frac{11}{24} u_{i+2} + \frac{3}{40} u_{i+3},
\end{equation} 
is introduced in \cite{Rathan} such that the finite difference numerical flux $\hat{g}_{i+1/2}^{\FD}$,
$$
   \hat{g}_{i+1/2}^{\FD} = \theta_0 \hat{g}^0_{i+1/2} + \theta_1 \hat{g}^1_{i+1/2} + \theta_2 \hat{g}^2_{i+1/2} + \theta_\C \hat{g}^\C_{i+1/2},
$$
is a convex combination of the numerical fluxes $\hat{g}^k_{i+1/2}$ \eqref{eq:numerical_flux_substencil} and $\hat{g}^\C_{i+1/2}$ \eqref{eq:numerical_flux_central} with 
$\theta_0 = \theta_2 = \frac{1}{6}, \, \theta_1 = \theta_\C = \frac{1}{3}$.
Note that the positive linear weights $\theta_k, \, k \in \{ 0, 1, 2, \text{C} \}$ can be set artificially as long as their sum is $1$ and $\theta_0 = \theta_2$.
With the central smoothness indicator
\begin{align*}
 \beta_\C = & \frac{4273}{20160} \left( u_{i-2} - 5 u_{i-1} + 10 u_i - 10 u_{i+1} + 5 u_{i+2} - u_{i+3} \right)^2 + \\
            & \frac{29}{345600} \left( 5 u_{i-2} + 11 u_{i-1} - 70 u_i + 94 u_{i+1} - 47 u_{i+2} + 7 u_{i+3} \right)^2 + \\
            & \frac{1}{3600} \left( 35 u_{i-2} - 139 u_{i-1} + 230 u_i - 206 u_{i+1} + 103 u_{i+2} - 23 u_{i+3} \right)^2 + \\
            & \frac{1}{576} \left( 7 u_{i-2} - 51 u_{i-1} + 134 u_i - 166 u_{i+1} + 99 u_{i+2} - 23 u_{i+3} \right)^2 + \\           
            & \frac{1}{2304} \left( 7 u_{i-2} - 56 u_{i-1} + 106 u_i - 76 u_{i+1} + 23 u_{i+2} - 4 u_{i+3} \right)^2 + \\
            & \frac{1}{9216} \left( 65 u_{i-2} - 353 u_{i-1} + 690 u_i - 602 u_{i+1} + 221 u_{i+2} - 21 u_{i+3} \right)^2 + \\
            & \frac{1}{9216} \left( 23 u_{i-2} - 63 u_{i-1} - 34 u_i + 186 u_{i+1} - 133 u_{i+2} + 21 u_{i+3} \right)^2 + \\
            & \frac{1}{2304} \left( 13 u_{i-2} - 28 u_{i-1} + 30 u_i - 28 u_{i+1} + 13 u_{i+2} \right)^2 + \\
            & \frac{2}{15} \left( u_{i-2} - 4 u_{i-1} + 6 u_i - 4 u_{i+1} + u_{i+2} \right)^2 + \\      
            & \frac{1}{1152} \left( u_{i-2} - 12 u_{i-1} + 22 u_i - 12 u_{i+1} + u_{i+2} \right)^2,
\end{align*} 
the new global smoothness indicator $\tau_6$ over the stencil $S^6$ is set as
$$
   \tau_6 = \left| \beta_\C - \frac{1}{24} ( 5 \beta_0 + 14 \beta_1 + 5 \beta_2) \right|.
$$
The nonlinear weights are defined by
\begin{equation} \label{eq:weights_CWENO}
 \omega^{\CWENO}_k = \frac{\alpha_k}{\sum_l \alpha_l}, \quad
 \alpha_k = \theta_k \left( 1 + \frac{\tau_6}{\beta_k + \epsilon} \right), \quad
 k, \, l \in \{ 0, 1, 2, \text{C} \}.
\end{equation}
Then the CWENO numerical flux is given by the following 
$$
   \hat{g}_{i+1/2} = \sum_{k \in \{ 0, 1, 2, \text{C} \}} \omega^{\CWENO}_k \hat{g}^k_{i+1/2}.
$$

Therefore, we approximate \eqref{eq:rd_discrete} by the following WENO scheme
$$
 \frac{du_{i}(t)}{dt} = \frac{\hat{g}_{i+1/2} - \hat{g}_{i-1/2}}{\Delta x^2} + R(u_{i}(t)),
$$
where $u_i(t)$ is the numerical approximation to $u(x_i,t)$.

\section{Numerical studies on convergence and wave speed} \label{sec:nr}
This section first show the accuracy of the centered finite difference (FD) and WENO methods for one-dimensional equations in terms of $L_1, \, L_2$ and $L_{\infty}$ error norms at the final time $t = T$.
We then compare the numerical results for large time steps and discuss the speed issue in the Newell–Whitehead–Segel equation.
For time discretization, we use the explicit third-order total variation diminishing Runge-Kutta method \cite{ShuOsherI}.
We follow the CFL condition in \cite{Liu} to set $\cfl = 0.4$ unless otherwise stated.
So the time step is $\Delta t = \cfl \cdot \Delta x^2$.
We choose $\epsilon = 10^{-6}$ for WENO-LSZ \cite{Liu}, $\epsilon = 10^{-30}$ for MWENO \cite{Hajipour} and $\epsilon = 10^{-40}$ for CWENO \cite{Rathan}.

\subsection{Convergence analysis} \label{sec:ca}
Numerical examples are presented to demonstrate the performance of FD and WENO methods explained in Section \ref{sec:weno}.
The numerical solutions are compared with the analytical solutions provided in Section \ref{sec:1d_rd}.

\begin{example} \label{ex:fisher}
We examine the accuracy for the Fisher's equation ($D=1, \, \rho=10^4$) on the domain $[-1, 5]$ with the initial condition given by \eqref{eq:fisher_exact} at $t=0$,
$$ u(x,0)= \frac{1}{\left[ 1 + \exp \left( \frac{100}{\sqrt{6}} x \right) \right]^2}, $$ 
and the Dirichlet boundary conditions $u(-1,t) = 1$ and $u(5,t) = 0$. 
The numerical solution is computed up to the final time $T=0.02$. 
The $L_1, \, L_2$ and $L_{\infty}$ errors versus $N$, are displayed in Tables \ref{tab:fisher_L1}, \ref{tab:fisher_L2} and \ref{tab:fisher_Linf}, respectively. 
It is clear that each scheme yield the convergent behavior as $N$ increases.
\end{example}

\begin{table}[htbp]
\renewcommand{\arraystretch}{1.1}
\scriptsize
\centering
\caption{$L_1$ error for Fisher’s equation.}      
\begin{tabular}{ccccc} 
\hline  
N & \multicolumn{4}{c}{$L_1$ error} \\ 
    \cline{2-5}
  & FD & WENO-LSZ & MWENO & CWENO \\
\hline
1200 & 1.072318e-4  & 1.073403e-4  & 1.091265e-4  & 1.072055e-4  \\  
2400 & 1.853247e-6  & 1.853247e-6  & 1.853624e-6  & 1.853175e-6  \\  
4800 & 2.970083e-8  & 2.970083e-8  & 2.970091e-8  & 2.970075e-8  \\
9600 & 4.670538e-10 & 4.670535e-10 & 4.670533e-10 & 4.670672e-10 \\ 
\hline
\end{tabular}
\label{tab:fisher_L1}
\end{table}

\begin{table}[htbp]
\renewcommand{\arraystretch}{1.1}
\scriptsize
\centering
\caption{$L_2$ error for Fisher’s equation.}      
\begin{tabular}{ccccc} 
\hline  
N & \multicolumn{4}{c}{$L_2$ error} \\ 
    \cline{2-5}
  & FD & WENO-LSZ & MWENO & CWENO \\
\hline
1200 & 7.512542e-4 & 7.512535e-4 & 7.645210e-4 & 7.510699e-4 \\  
2400 & 1.298020e-5 & 1.298020e-5 & 1.298284e-5 & 1.297970e-5 \\  
4800 & 2.080026e-7 & 2.080026e-7 & 2.080031e-7 & 2.080020e-7 \\
9600 & 3.270464e-9 & 3.270464e-9 & 3.270464e-9 & 3.270464e-9 \\ 
\hline
\end{tabular}
\label{tab:fisher_L2}
\end{table}

\begin{table}[htbp]
\renewcommand{\arraystretch}{1.1}
\scriptsize
\centering
\caption{$L_{\infty}$ error for Fisher’s equation.}      
\begin{tabular}{ccccc} 
\hline
N & \multicolumn{4}{c}{$L_{\infty}$ error} \\ 
    \cline{2-5}
  & FD & WENO-LSZ & MWENO & CWENO \\
\hline
1200 & 7.795743e-3 & 7.795739e-3 & 7.933331e-3 & 7.793864e-3 \\  
2400 & 1.346346e-4 & 1.346346e-4 & 1.346616e-4 & 1.346296e-4 \\ 
4800 & 2.157525e-6 & 2.157525e-6 & 2.157527e-6 & 2.157521e-6 \\
9600 & 3.392635e-8 & 3.392635e-8 & 3.392633e-8 & 3.392636e-8 \\  
\hline
\end{tabular}
\label{tab:fisher_Linf}
\end{table}

\begin{example} \label{ex:zeldovich}
We consider the Zeldovich equation ($D=1, \, \rho=9000$) on the domain $[-1, 5]$.
The initial condition \eqref{eq:zeldovich_exact} is 
$$ u(x,0)= \frac{1}{1 + \exp \left( 30 \sqrt{5} x  \right)}, $$
and the boundary conditions $u(-1,t) = 1$ and $u(5,t) = 0$. 
The final time $T=0.06$.
The $L_1, \, L_2$ and $L_{\infty}$ errors versus $N$, are displayed in Tables \ref{tab:zeldovich_L1}, \ref{tab:zeldovich_L2} and \ref{tab:zeldovich_Linf}, respectively.
We observe that all schemes converge with the increase of $N$, but the errors of the CWENO scheme are larger than FD, WENO-LSZ and MWENO.
\end{example}

\begin{table}[htbp]
\renewcommand{\arraystretch}{1.1}
\scriptsize
\centering
\caption{$L_1$ error for Zeldovich equation.}      
\begin{tabular}{ccccc} 
\hline  
N & \multicolumn{4}{c}{$L_1$ error} \\ 
    \cline{2-5} 
  & FD & WENO-LSZ & MWENO & CWENO \\
\hline
1200 & 7.722126e-7  & 9.088312e-7  & 3.012620e-7  & 1.655857e-6  \\  
2400 & 1.224838e-8  & 1.234412e-8  & 1.139223e-8  & 1.469268e-8  \\  
4800 & 1.909269e-10 & 1.909286e-10 & 1.901880e-10 & 1.966083e-10 \\
9600 & 6.038550e-12 & 6.038878e-12 & 6.037716e-12 & 6.071350e-12 \\ 
\hline
\end{tabular}
\label{tab:zeldovich_L1}
\end{table}

\begin{table}[htbp]
\renewcommand{\arraystretch}{1.1}
\scriptsize
\centering
\caption{$L_2$ error for Zeldovich equation.}      
\begin{tabular}{ccccc} 
\hline  
N & \multicolumn{4}{c}{$L_2$ error} \\ 
    \cline{2-5}
  & FD & WENO-LSZ & MWENO & CWENO \\
\hline
1200 & 6.357072e-6  & 7.477862e-6  & 2.492258e-6  & 1.360983e-5  \\  
2400 & 1.008014e-7  & 1.015865e-7  & 9.377596e-8  & 1.208601e-7  \\  
4800 & 1.571049e-9  & 1.571063e-9  & 1.564984e-9  & 1.617638e-9  \\
9600 & 4.926294e-11 & 4.926314e-11 & 4.925981e-11 & 4.936196e-11 \\ 
\hline
\end{tabular}
\label{tab:zeldovich_L2}
\end{table}

\begin{table}[htbp]
\renewcommand{\arraystretch}{1.1}
\scriptsize
\centering
\caption{$L_{\infty}$ error for Zeldovich equation.}      
\begin{tabular}{ccccc} 
\hline
N & \multicolumn{4}{c}{$L_{\infty}$ error} \\ 
    \cline{2-5}  
  & FD & WENO-LSZ & MWENO & CWENO \\
\hline
1200 & 7.902828e-5  & 9.281344e-5  & 3.144861e-5  & 1.684816e-4  \\  
2400 & 1.252654e-6  & 1.262307e-6  & 1.166094e-6  & 1.500097e-6  \\ 
4800 & 1.954562e-8  & 1.954579e-8  & 1.947092e-8  & 2.012174e-8  \\
9600 & 6.088552e-10 & 6.088578e-10 & 6.088167e-10 & 6.100815e-10 \\  
\hline
\end{tabular}
\label{tab:zeldovich_Linf}
\end{table}

\begin{example} \label{ex:nws}
We solve the Newell–Whitehead–Segel equation ($D=1, \, \rho=5000, \, \alpha=2$) on the domain $[-1, 5]$ with the initial condition \eqref{eq:nws_exact}
$$ u(x,0) = -\frac{1}{2} \tanh \left( 25 x \right) + \frac{1}{2}, $$ 
and the boundary conditions $u(-1,t) = 1$ and $u(5,t) = 0$. 
The final time $T=0.028$.
The $L_1, \, L_2$ and $L_{\infty}$ errors versus $N$, are displayed in Tables \ref{tab:nws_L1}, \ref{tab:nws_L2} and \ref{tab:nws_Linf}, respectively.
Note that none of the schemes is convergent because of the speed issue, which will be investigated in detail in the following section.
\end{example}

\begin{table}[htbp]
\renewcommand{\arraystretch}{1.1}
\scriptsize
\centering
\caption{$L_1$ error for Newell–Whitehead–Segel equation ($\alpha=2$).}      
\begin{tabular}{ccccc} 
\hline  
N & \multicolumn{4}{c}{$L_1$ error} \\
    \cline{2-5}
  & FD & WENO-LSZ & MWENO & CWENO \\
\hline
1200 & 0.016920 & 0.032764 & 0.016920 & 0.016916 \\  
2400 & 0.017003 & 0.017003 & 0.017003 & 0.017002 \\  
4800 & 0.016991 & 0.016991 & 0.016991 & 0.016991 \\
9600 & 0.016977 & 0.016977 & 0.016977 & 0.016977 \\ 
\hline
\end{tabular}
\label{tab:nws_L1}
\end{table}

\begin{table}[htbp]
\renewcommand{\arraystretch}{1.1}
\scriptsize
\centering
\caption{$L_2$ error for Newell–Whitehead–Segel equation ($\alpha=2$).}      
\begin{tabular}{ccccc} 
\hline  
N & \multicolumn{4}{c}{$L_2$ error} \\ 
    \cline{2-5}
  & FD & WENO-LSZ & MWENO & CWENO \\
\hline
1200 & 0.103038 & 0.109786 & 0.103039 & 0.103018 \\  
2400 & 0.103411 & 0.103411 & 0.103411 & 0.103406 \\  
4800 & 0.103350 & 0.103350 & 0.103350 & 0.103349 \\
9600 & 0.103281 & 0.103281 & 0.103281 & 0.103280 \\ 
\hline
\end{tabular}
\label{tab:nws_L2}
\end{table}

\begin{table}[htbp]
\renewcommand{\arraystretch}{1.1}
\scriptsize
\centering
\caption{$L_{\infty}$ error for Newell–Whitehead–Segel equation ($\alpha=2$).}      
\begin{tabular}{ccccc} 
\hline
N & \multicolumn{4}{c}{$L_{\infty}$ error} \\ 
    \cline{2-5}
  & FD & WENO-LSZ & MWENO & CWENO \\
\hline
1200 & 0.864944 & 0.864944 & 0.864946 & 0.864862 \\  
2400 & 0.866544 & 0.866544 & 0.866544 & 0.866524 \\ 
4800 & 0.866379 & 0.866379 & 0.866379 & 0.866373 \\
9600 & 0.866059 & 0.866059 & 0.866059 & 0.866057 \\  
\hline
\end{tabular}
\label{tab:nws_Linf}
\end{table}

\begin{example} \label{ex:bistable}
We test the accuracy for the bistable equation ($D=1, \, \rho=10^4, \, \beta=0.2$). 
The initial condition \eqref{eq:bistable_exact} is given by
$$ u(x,0) = \frac{3}{5} + \frac{2}{5} \tanh \left( 20 \sqrt{2} x \right). $$ 
The computational domain is $[-5, 1]$ with the boundary conditions $u(-5,t) = 0.2$ and $u(1,t) = 1$, and the final time is $T=0.05$.
The $L_1, \, L_2$ and $L_{\infty}$ errors versus $N$, are displayed in Tables \ref{tab:bistable_L1}, \ref{tab:bistable_L2} and \ref{tab:bistable_Linf}, respectively.
The numerical simulations match the exact solution well up to at least $4$ decimal digits of accuracy, even though the errors begin to increase after $N=2400$, which is likely due to the effect of the large CFL number in time integration.    
\end{example} 

\begin{table}[htbp]
\renewcommand{\arraystretch}{1.1}
\scriptsize
\centering
\caption{$L_1$ error for bistable equation.}      
\begin{tabular}{ccccc} 
\hline  
N & \multicolumn{4}{c}{$L_1$ error} \\
    \cline{2-5}
  & FD & WENO-LSZ & MWENO & CWENO \\
\hline
1200 & 8.293439e-7 & 8.394313e-7 & 8.005169e-7 & 8.941725e-7 \\  
2400 & 3.067402e-8 & 3.067864e-8 & 3.064174e-8 & 3.050023e-8 \\  
4800 & 3.592002e-8 & 3.591735e-8 & 3.591643e-8 & 3.564344e-8 \\
9600 & 7.229070e-8 & 7.228750e-8 & 7.256824e-8 & 7.311459e-8 \\ 
\hline
\end{tabular}
\label{tab:bistable_L1}
\end{table}

\begin{table}[htbp]
\renewcommand{\arraystretch}{1.1}
\scriptsize
\centering
\caption{$L_2$ error for bistable equation.}      
\begin{tabular}{ccccc} 
\hline  
N & \multicolumn{4}{c}{$L_2$ error} \\ 
    \cline{2-5}
  & FD & WENO-LSZ & MWENO & CWENO \\
\hline
1200 & 6.247502e-6 & 6.323723e-6 & 6.028918e-6 & 6.738280e-6 \\  
2400 & 2.302114e-7 & 2.302445e-7 & 2.299654e-7 & 2.289235e-7 \\  
4800 & 2.689060e-7 & 2.688858e-7 & 2.688790e-7 & 2.668354e-7 \\
9600 & 5.412034e-7 & 5.411788e-7 & 5.432810e-7 & 5.473708e-7 \\ 
\hline
\end{tabular}
\label{tab:bistable_L2}
\end{table}

\begin{table}[htbp]
\renewcommand{\arraystretch}{1.1}
\scriptsize
\centering
\caption{$L_{\infty}$ error for bistable equation.}      
\begin{tabular}{clcrlcrlcrlc} 
\hline
N & \multicolumn{4}{c}{$L_{\infty}$ error} \\ 
    \cline{2-5}
  & FD & WENO-LSZ & MWENO & CWENO \\
\hline 
1200 & 7.072738e-5 & 7.158588e-5 & 6.825260e-5 & 7.635299e-5 \\  
2400 & 2.601050e-6 & 2.601417e-6 & 2.598225e-6 & 2.586766e-6 \\ 
4800 & 3.028937e-6 & 3.028708e-6 & 3.028633e-6 & 3.005615e-6 \\
9600 & 6.095777e-6 & 6.095498e-6 & 6.119178e-6 & 6.165241e-6 \\  
\hline
\end{tabular}
\label{tab:bistable_Linf}
\end{table}

\begin{example} \label{ex:lotka_volterra}
We conclude this subsection with the Lotka-Volterra competition-diffusion system \eqref{eq:lotka_volterra} ($D=1, \, \rho=7000$) with the initial condition \eqref{eq:lotka_volterra_exact},
\begin{align*}
 u(x,0) &= \frac{1}{2} \left[ 1 + \tanh \left( 5 \sqrt{105} x \right) \right], \\
 v(x,0) &= \frac{3}{4} \left[ 1 - \tanh \left( 5 \sqrt{105} x \right) \right]^2.
\end{align*}
The computational domain is $[-1, 5]$, and the boundary conditions are $(u,v)(-1,t) = (0,3)$ and $(u,v)(5,t) = (1,0)$. 
We compute the numerical solution at the final time $T=0.1$. 
The $L_1, \, L_2$ and $L_{\infty}$ errors versus $N$, are displayed in Tables \ref{tab:lotka_volterra_L1}, \ref{tab:lotka_volterra_L2} and \ref{tab:lotka_volterra_Linf}, respectively, which illustrates that all schemes exhibit convergence to the exact solutions.    
\end{example}

\begin{table}[htbp]
\renewcommand{\arraystretch}{1.1}
\scriptsize
\centering
\caption{$L_1$ error for Lotka-Volterra competition-diffusion system.}      
\begin{tabular}{cccccc} 
\hline  
& N & \multicolumn{4}{c}{$L_1$ error} \\ 
      \cline{3-6}
&   & FD & WENO-LSZ & MWENO & CWENO \\
\hline
\multirow{3}{*}{$u$} & 1500  & 1.666967e-6  & 7.659615e-6  & 2.109070e-5  & 1.738599e-5  \\  
                     & 3000  & 2.695194e-8  & 4.032780e-8  & 5.555932e-8  & 8.358504e-8  \\  
                     & 6000  & 4.224317e-10 & 4.298484e-10 & 6.095723e-11 & 6.747933e-10 \\
\hline
\multirow{3}{*}{$v$} & 1500  & 5.013286e-6 & 2.299186e-5 & 6.328052e-5  & 5.21364357e-5 \\  
                     & 3000  & 8.105937e-8 & 1.211846e-7 & 1.664837e-7  & 2.50861548e-7 \\  
                     & 6000  & 1.270529e-9 & 1.292779e-9 & 1.860878e-10 & 2.02741127e-9 \\
\hline
\end{tabular}
\label{tab:lotka_volterra_L1}
\end{table}

\begin{table}[htbp]
\renewcommand{\arraystretch}{1.1}
\scriptsize
\centering
\caption{$L_2$ error for Lotka-Volterra competition-diffusion system.}      
\begin{tabular}{cccccc} 
\hline  
& N & \multicolumn{4}{c}{$L_2$ error} \\ 
      \cline{3-6}
&   & FD & WENO-LSZ & MWENO & CWENO \\
\hline
\multirow{3}{*}{$u$} & 1500 & 1.690017e-5 & 7.763889e-5 & 2.137994e-4  & 1.761342e-4 \\  
                     & 3000 & 2.732035e-7 & 4.086730e-7 & 5.626008e-7  & 8.467678e-7 \\  
                     & 6000 & 4.281759e-9 & 4.356861e-9 & 6.221991e-10 & 6.837337e-9 \\
\hline
\multirow{3}{*}{$v$} & 1500 & 5.568829e-5 & 2.551775e-4 & 7.021555e-4 & 5.785293e-4 \\ 
                     & 3000 & 9.003274e-7 & 1.345496e-6 & 1.846491e-6 & 2.784076e-6 \\
                     & 6000 & 1.411056e-8 & 1.435737e-8 & 2.089350e-9 & 2.250653e-8 \\
\hline
\end{tabular}
\label{tab:lotka_volterra_L2}
\end{table}

\begin{table}[htbp]
\renewcommand{\arraystretch}{1.1}
\scriptsize
\centering
\caption{$L_{\infty}$ error for Lotka-Volterra competition-diffusion system.}      
\begin{tabular}{cccccc} 
\hline
& N & \multicolumn{4}{c}{$L_{\infty}$ error} \\ 
      \cline{3-6}
&  & FD & WENO-LSZ & MWENO & CWENO \\
\hline
\multirow{3}{*}{$u$} & 1500 & 2.579368e-4 & 1.180368e-3 & 3.245438e-3 & 2.677192e-3 \\  
                     & 3000 & 4.170033e-6 & 6.227577e-6 & 8.524106e-6 & 1.288152e-5 \\ 
                     & 6000 & 6.549232e-8 & 6.663192e-8 & 1.004506e-8 & 1.042878e-7 \\
\hline
\multirow{3}{*}{$v$} & 1500 & 9.013902e-4 & 4.149239e-3 & 1.145263e-2 & 9.411484e-3 \\  
                     & 3000 & 1.469798e-5 & 2.203628e-5 & 3.059522e-5 & 4.572058e-5 \\ 
                     & 6000 & 2.302456e-7 & 2.343120e-7 & 3.534330e-8 & 3.687668e-7 \\ 
\hline
\end{tabular}
\label{tab:lotka_volterra_Linf}
\end{table}

\subsection{Time step: stability and convergence} \label{sec:timestepping}
It is well known that the explicit temporal method is stable only for very small time steps of magnitude $\Delta x^2$ as in Subsection \ref{sec:ca}, when there is a diffusion term involved.
One advantage of the CWENO scheme is that it is possible to choose a larger time step than FD, WENO-LSZ and MWENO if we choose a relatively large value of $\Delta x$.
In order to investigate the time stepping for the traveling wave solutions with sharp fronts, we consider the following numerical examples for different equations. 

\begin{figure}[htbp]
\centering
\includegraphics[width=0.325\textwidth]{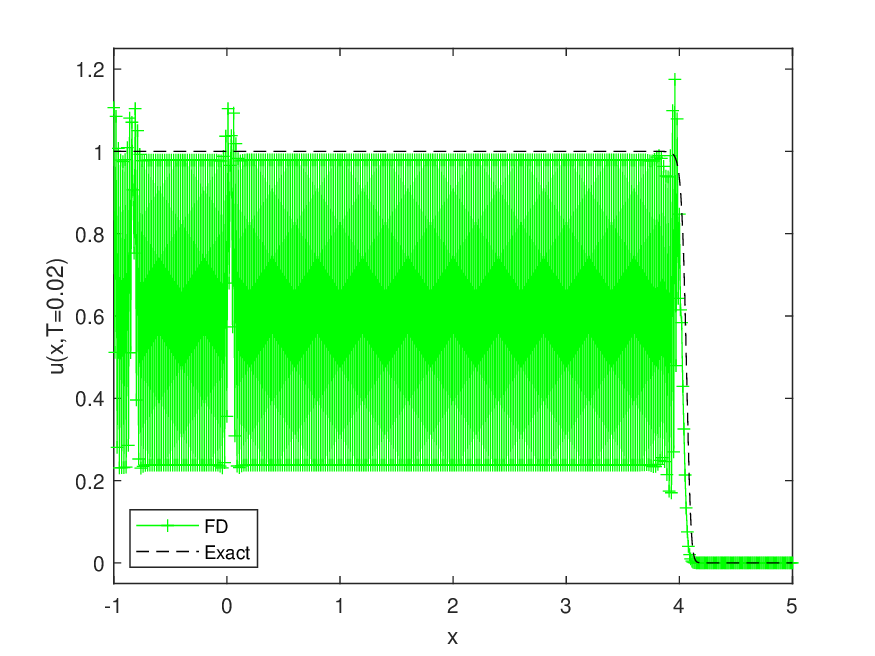}
\includegraphics[width=0.325\textwidth]{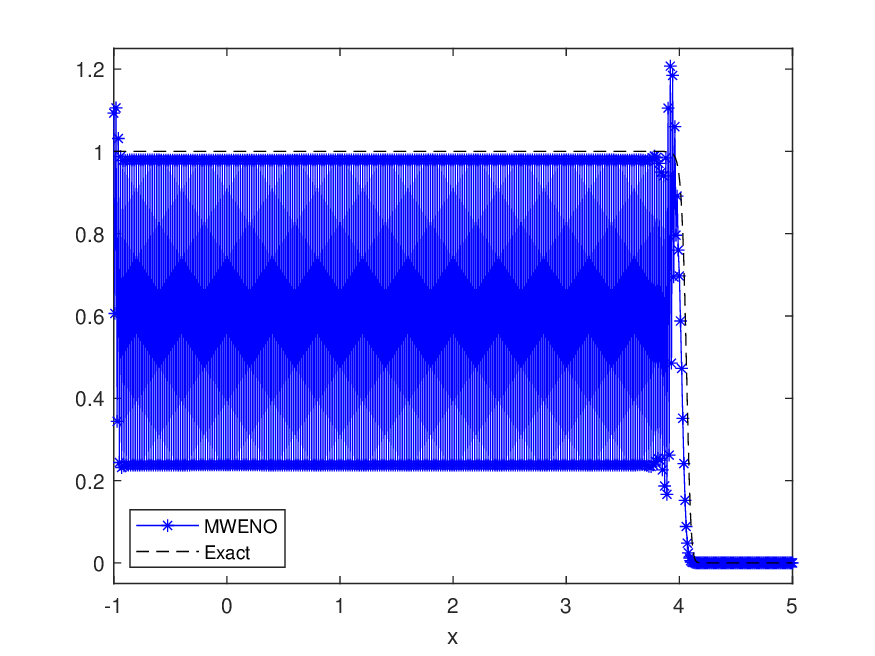}
\includegraphics[width=0.325\textwidth]{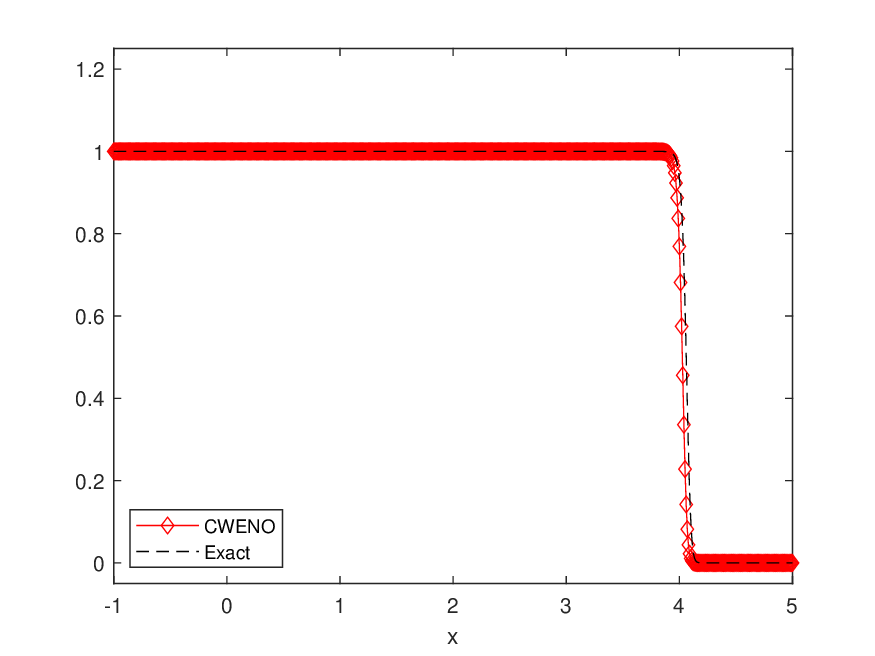}
\caption{Solution profiles for Fisher's equation at $T=0.02$ computed by FD (left), MWENO (middle) and CWENO (right) with $N=600$.
The dashed black line is the exact solution.}
\label{fig:fisher_time_step}
\end{figure}

For the Fisher's equation ($D = 1, \, \rho=10^4$) with the initial condition $u(x,0)$ in \eqref{eq:fisher_exact}, the computational domain $[-1, \, 5]$ is divided into $N=600$ uniform cells.
The numerical solution is computed up to the final time $T=0.02$.
The simulation run by WENO-LSZ blows up at $t=0.00248$.
Figure \ref{fig:fisher_time_step} shows that the solutions by FD, MWENO and CWENO are all stable, but only CWENO yields non-oscillatory and accurate solution.  

\begin{figure}[htbp]
\centering
\includegraphics[width=0.495\textwidth]{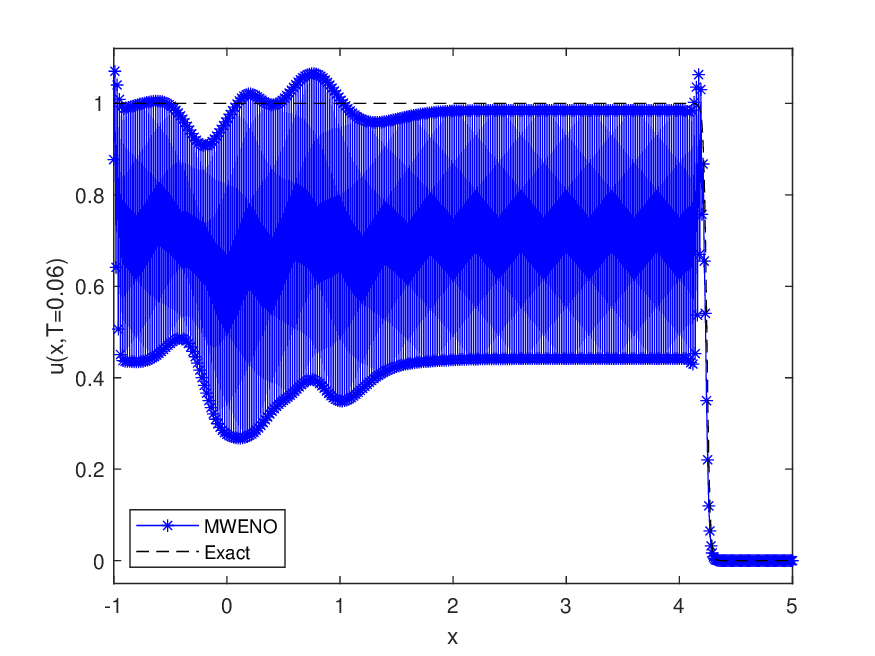}
\includegraphics[width=0.495\textwidth]{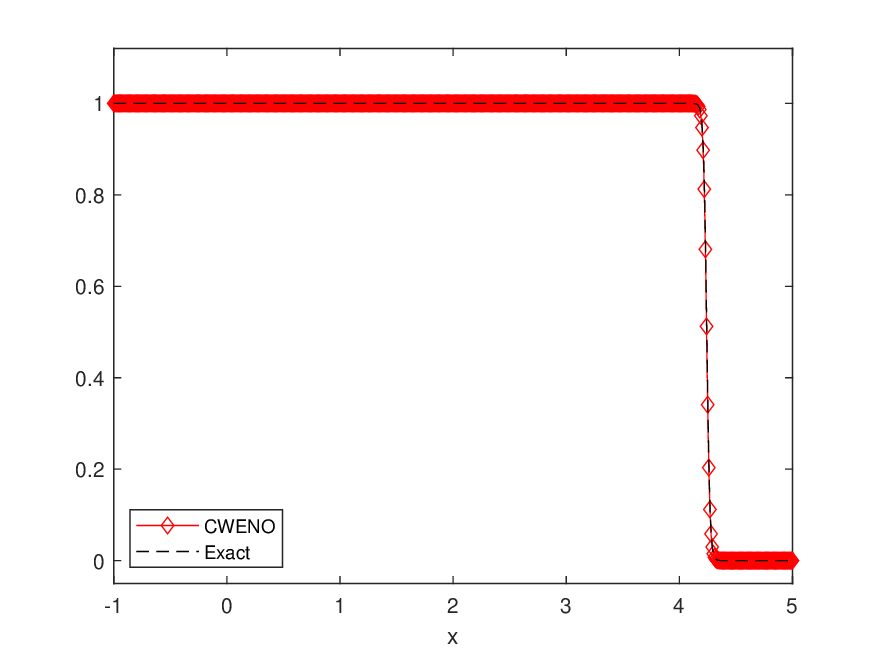}
\caption{Solution profiles for Zeldovich equation at $T=0.06$ computed by MWENO (left) and CWENO (right) with $N=600$.
The dashed black line is the exact solution.}
\label{fig:zeldovich_time_step}
\end{figure}

For the Zeldovich equation ($D = 1, \, \rho=10^4$) with the initial condition $u(x,0)$ in \eqref{eq:zeldovich_exact}, we divide the computational domain $[-1, \, 5]$ into $N=600$ uniform cells.
The final time is $T=0.06$. 
The numerical solutions by FD and WENO-LSZ blow up at the time $t=0.00184$ and $t=0.00172$, respectively.
However, the solutions by MWENO and CWENO are stable with the solution by CWENO non-oscillatory and accurate, as shown in Figure \ref{fig:zeldovich_time_step}.

\begin{figure}[htbp]
\centering
\includegraphics[width=0.325\textwidth]{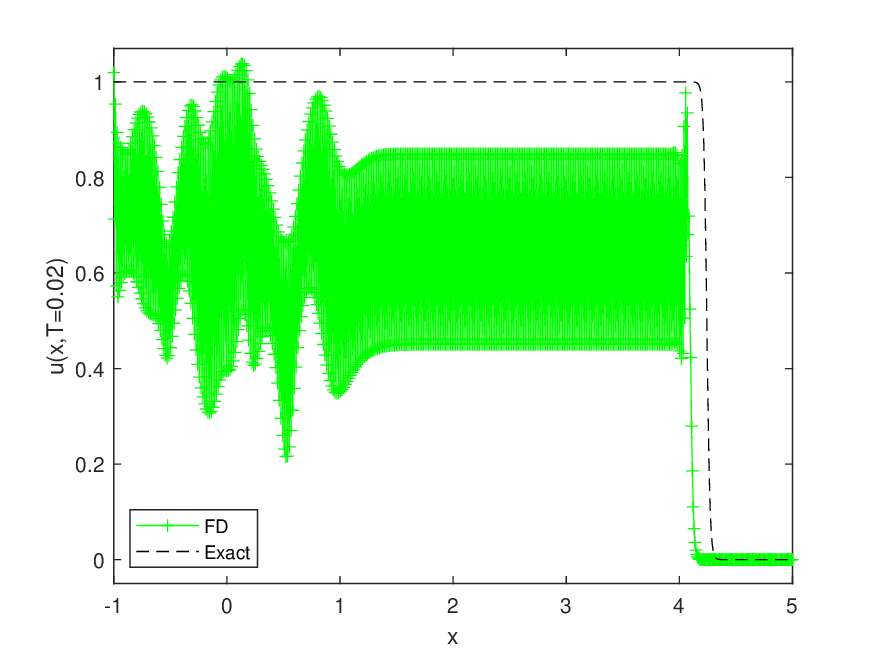}
\includegraphics[width=0.325\textwidth]{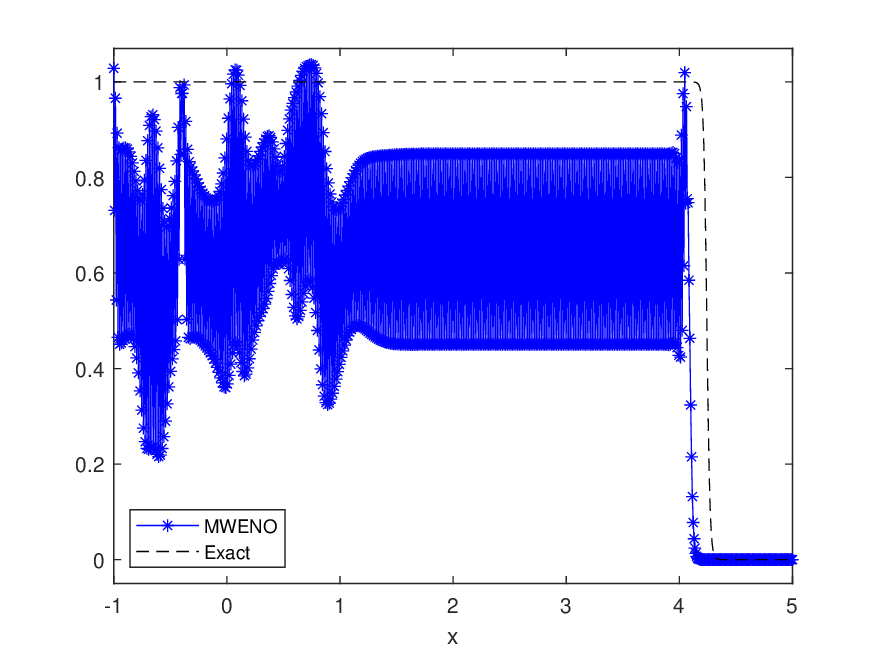}
\includegraphics[width=0.325\textwidth]{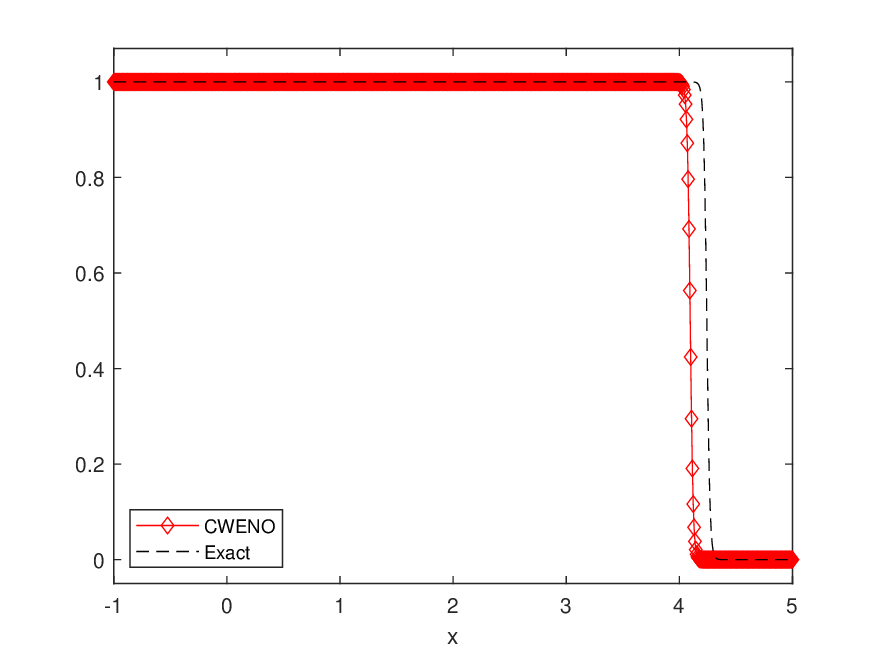}
\caption{Solution profiles for Newell-Whitehead-Segel equation at $T=0.02$ computed by FD (left), MWENO (middle) and CWENO (right) with $N=800$.
The dashed black line is the exact solution.}
\label{fig:nws_time_step}
\end{figure}

For the Newell-Whitehead-Segel equation ($D = 1, \, \rho=10^4, \, \alpha=2$), the initial condition $u(x,0)$ is given by \eqref{eq:nws_exact}, and the computational domain is $[-1, \, 5]$.
We take $N=800$ and run the numerical solution up to the final time $T=0.02$.
The numerical solution by WENO-LSZ blows up at $t=0.0008775$.
Figure \ref{fig:nws_time_step} plots the numerical solutions by FD, MWENO and CWENO, where only the solution by CWENO is non-oscillatory and accurate. 
The solutions with FD and MWENO are similarly inaccurate and oscillatory, even though they are stable.
For this case, we also observe that the phase error is noticeable. 
That is, the wave front of the numerical solution lags behind the exact wave front, indicating that the numerical speed of the wave front is slower than the exact one. 
This phenomenon, which is consistent with Example \ref{ex:nws}, is also observed with FD and MWENO, despite their high oscillatory behaviors.

\begin{figure}[htbp]
\centering
\includegraphics[width=0.325\textwidth]{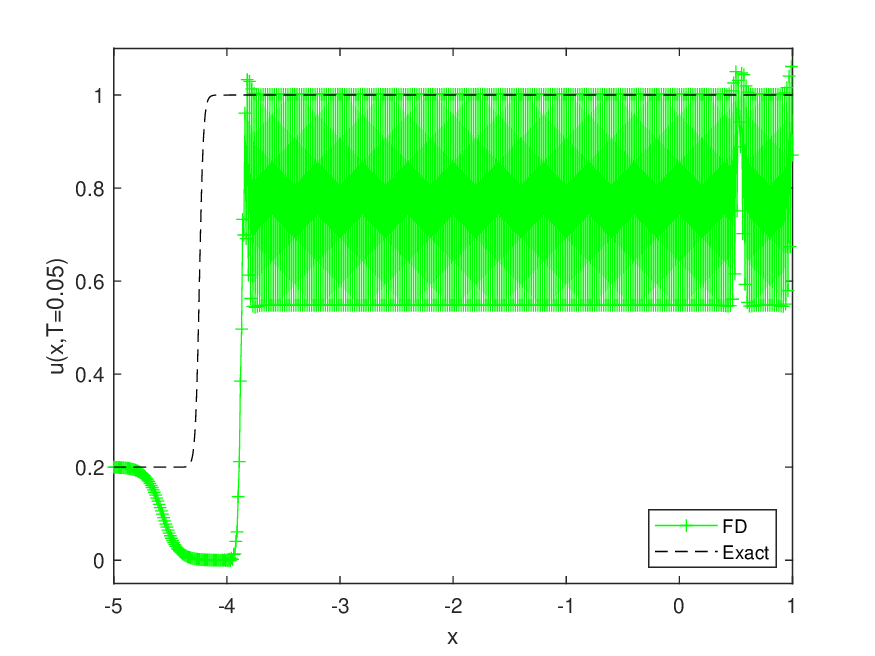}
\includegraphics[width=0.325\textwidth]{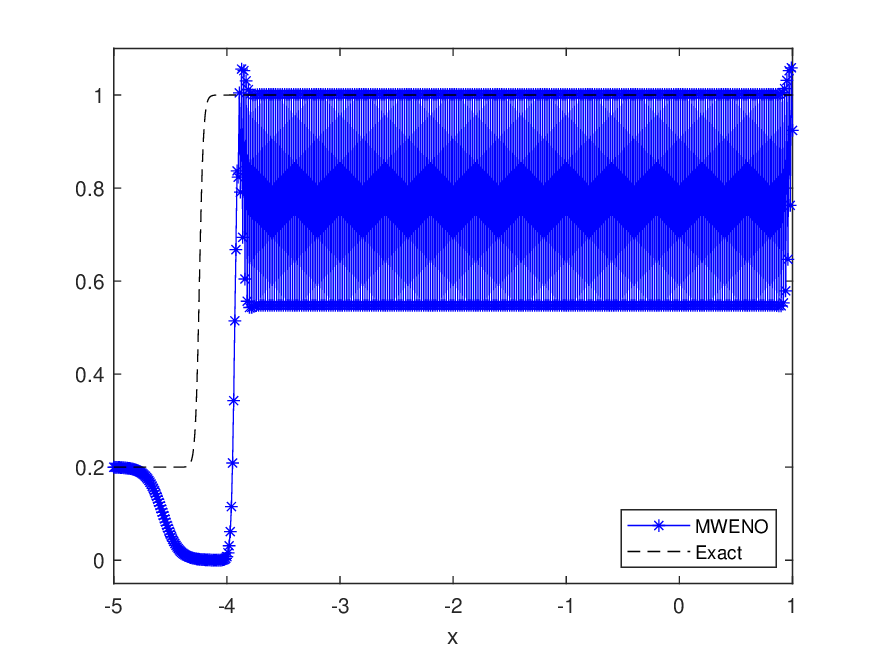}
\includegraphics[width=0.325\textwidth]{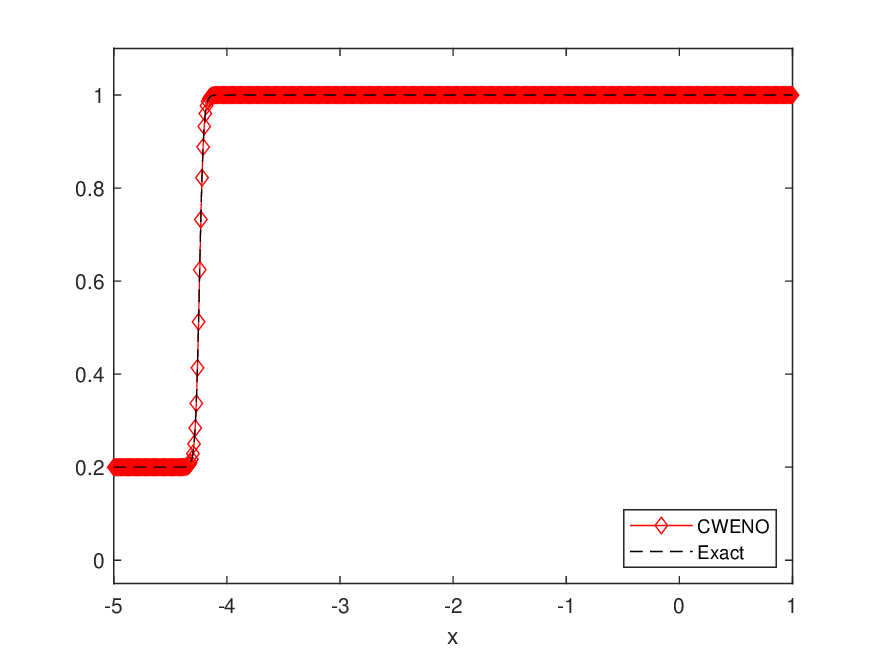}
\caption{Solution profiles for bistable equation at $T=0.02$ computed by FD (left), MWENO (middle) and CWENO (right) with $N=600$.
The dashed black line is the exact solution.}
\label{fig:bistable_time_step}
\end{figure}

For the bistable equation ($D = 1, \, \rho=10^4, \, \beta=0.2$) with the initial condition $u(x,0)$ in \eqref{eq:bistable_exact} for $t=0$, the computational domain is $[-5, \, 1]$ with $N=600$ uniform cells.
Figure \ref{fig:bistable_time_step} shows the numerical solutions by FD, MWENO and CWENO at the final time $T=0.02$.
Only the solution using CWENO is both stable and accurate. 
In contrast, the solutions obtained by FD and MWENO are stable but inaccurate, especially regarding the phase error of the traveling wave fronts.
Note that the numerical solution by WENO-LSZ blows up at $t=0.00276$.

\begin{figure}[htbp]
\centering
\includegraphics[width=0.325\textwidth]{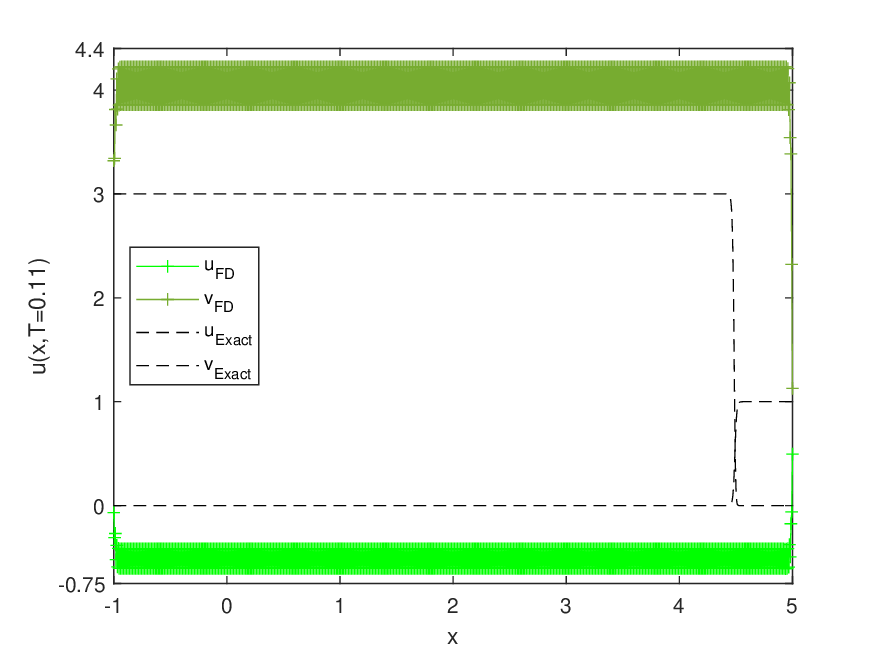}
\includegraphics[width=0.325\textwidth]{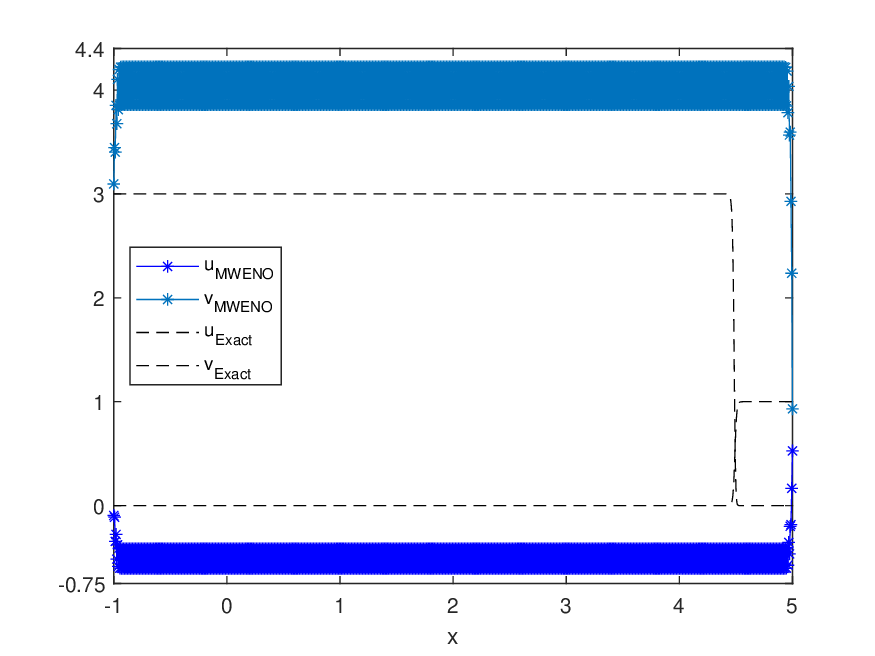}
\includegraphics[width=0.325\textwidth]{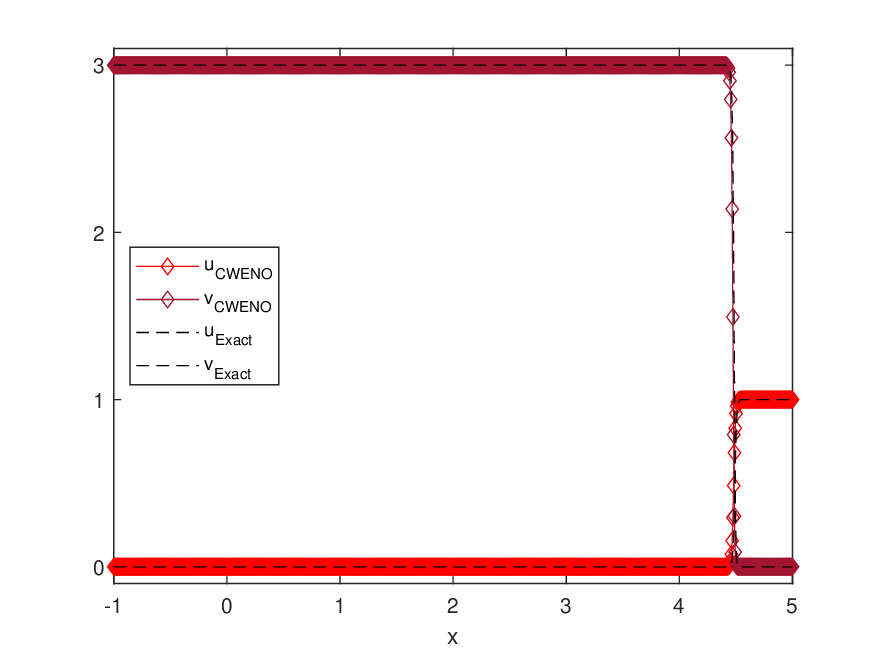}
\caption{Solution profiles for Lotka-Volterra competition-diffusion system at $T=0.11$ computed by FD (left), MWENO (middle) and CWENO (right) with $N=900$.
The dashed black line is the exact solution.}
\label{fig:lv_time_step}
\end{figure}

We consider the Lotka-Volterra competition-diffusion system ($D = 1, \, \rho=10^4$).
The initial condition $(u,v)(x,0)$ is given by \eqref{eq:lotka_volterra_exact} for $t=0$,
The computational domain $[-1, \, 5]$ is discretized with $N=900$ uniform cells and the final time is $T=0.11$.
We present the numerical results approximated by FD, MWENO and CWENO at the final time, as shown in Figure \ref{fig:lv_time_step}, where only the solution by CWENO is both stable and accurate.
Again the numerical solution by WENO-LSZ blows up at the time $t=0.00082$.

Table \ref{tab:comparison} provides the summary of the performance of each method. 
As shown in the table, only CWENO method yields stable and accurate solutions for the considered equations.

\begin{table}[htbp]
\centering
\begin{tabular}{|c|c|c|c|c|}
\hline
       & \textbf{FD} & \textbf{WENO-LSZ} & \textbf{MWENO} & \textbf{CWENO} \\
\hline
Fisher & Stable/Inaccurate & Unstable & Stable/Inaccurate & Stable/Accurate \\
\hline
Zeldovich & Unstable & Unstable & Stable/Inaccurate & Stable/Accurate \\
\hline
Newell & Stable/Inaccurate & Unstable & Stable/Inaccurate & Stable/Accurate \\
\hline
Bistable & Stable/Inaccurate & Unstable & Stable/Inaccurate & Stable/Accurate \\
\hline
Lotka-Volterra & Stable/Inaccurate & Unstable & Stable/Inaccurate & Stable/Accurate \\
\hline
\end{tabular}
\caption{Stability and accuracy: comparison of methods with different equations: FD, WENO-LSZ, MWENO, and CWENO methods. 
The table shows that only CWENO method yields stable and accurate solutions for the considered equations.}
\label{tab:comparison}
\end{table}

\subsection{Convergence of numerical speed} \label{sec:speeding}
In this subsection, we discuss the speed issue related to the stability and convergence of the methods. 
Specifically, we will demonstrate that for the Newell-Whitehead-Segel equation, arbitrarily decreasing the CFL number with a given $N$ does not ensure that the numerical speed of the wave front converges to the exact speed.
In Figure \ref{fig:nws_time_step}, the sharp front computed with each method clearly travels at the wrong speed, resulting in a lag behind the correct location. 
In order to investigate how the lag behaves with different values of $N$ and $\alpha$, we first refine the grid and increase $N$ for different values of $\alpha$ to evaluate the performance of CWENO.

\begin{figure}[htbp]
\centering
\includegraphics[width=0.325\textwidth]{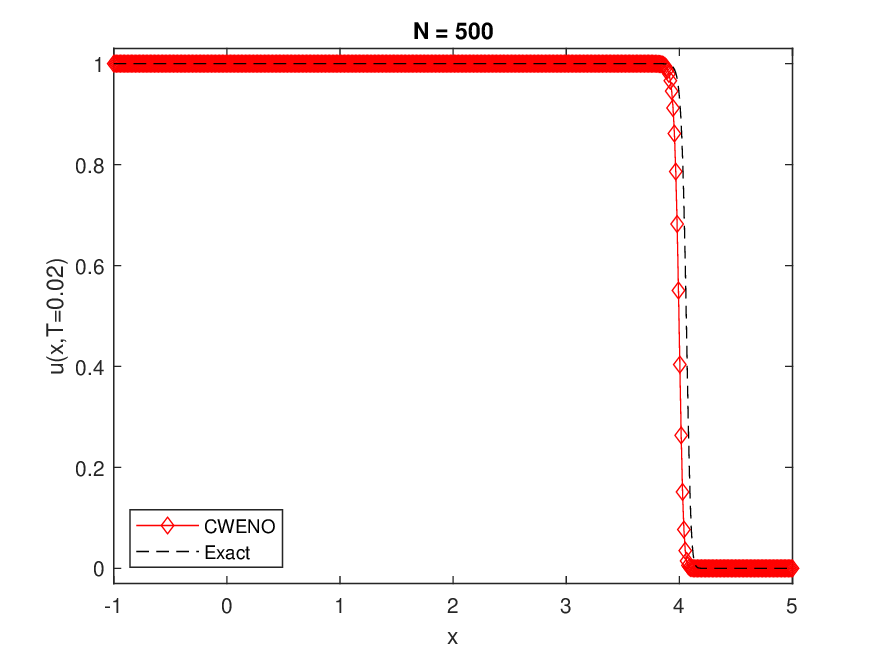}
\includegraphics[width=0.325\textwidth]{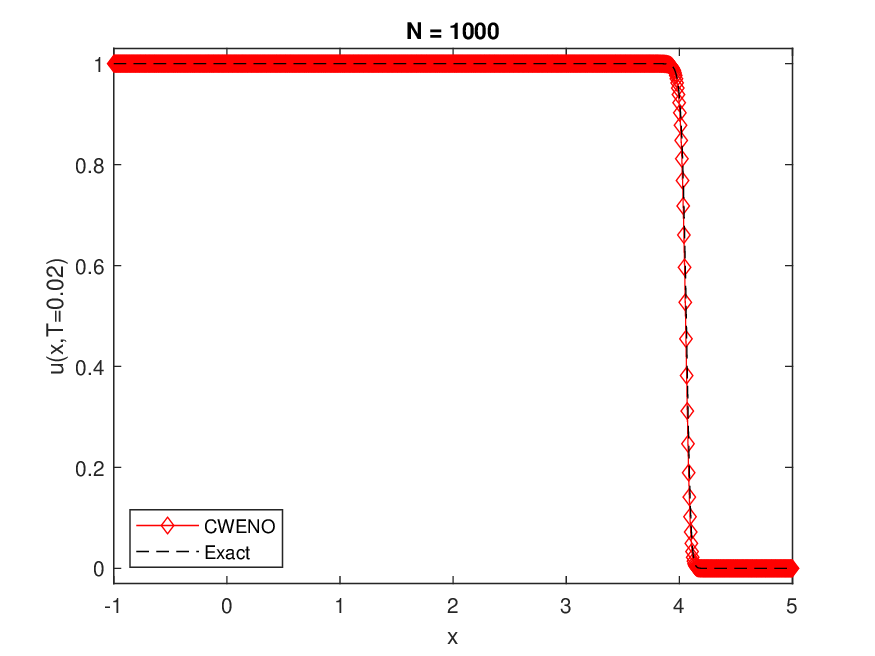}
\includegraphics[width=0.325\textwidth]{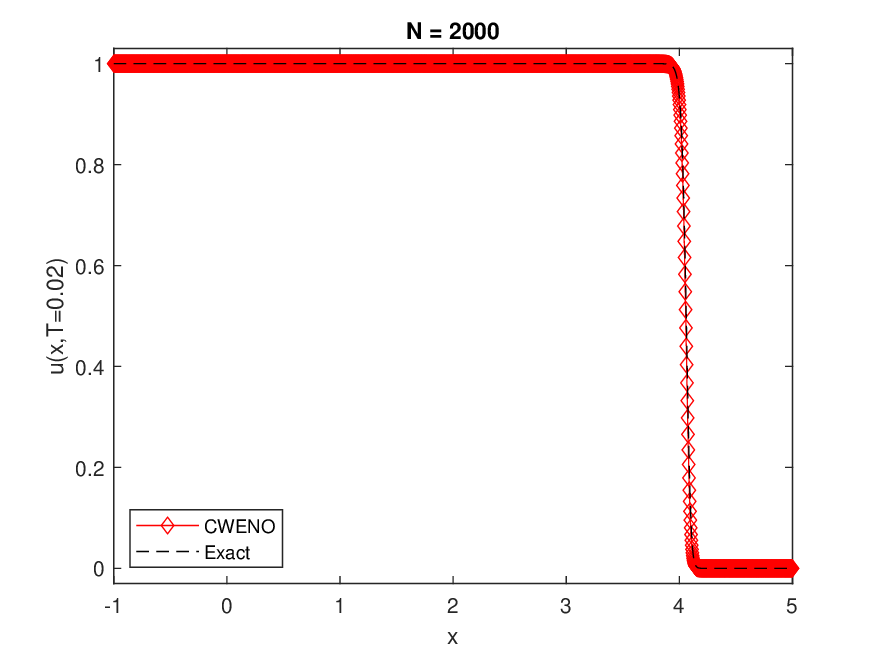}
\caption{Solution profiles for Newell-Whitehead-Segel equation ($\alpha=1$) at $T=0.02$ computed by CWENO with $N=500$ (left), $1000$ (middle) and $2000$ (right).
The dashed black line is the exact solution.}
\label{fig:nws_alpha1_N}
\end{figure}

Figure \ref{fig:nws_alpha1_N} shows solution profiles for the Newell-Whitehead-Segel equation ($D = 1, \, \rho=10^4$) with $\alpha = 1$ (or equivalently, Fisher's equation) at $T = 0.02$ computed by CWENO. 
As shown in the figure, we see the convergence to the exact solution as $N$ increases and the error becomes smaller as shown in Table \ref{tab:nws_error_N}.

\begin{figure}[htbp]
\centering
\includegraphics[width=0.325\textwidth]{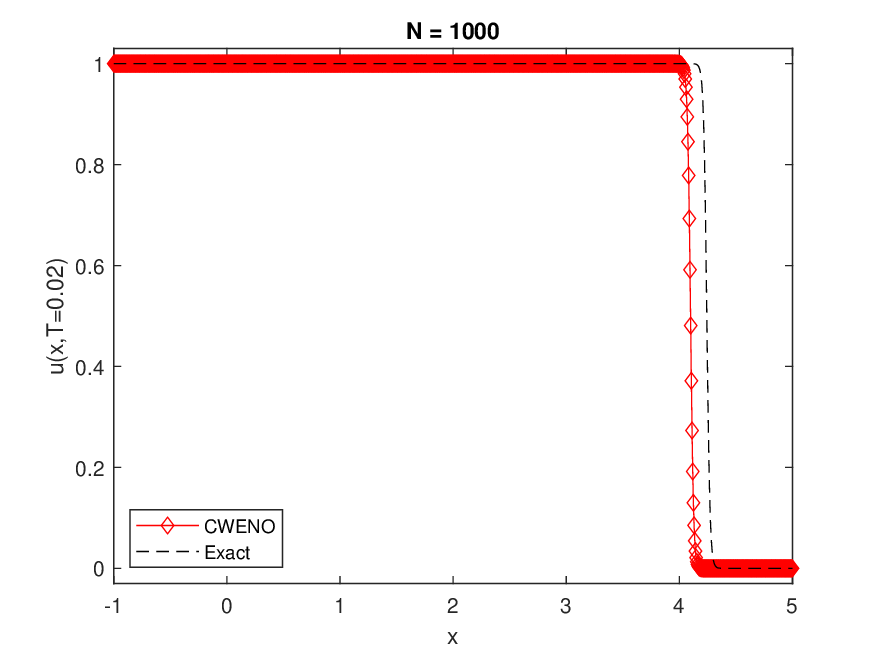}
\includegraphics[width=0.325\textwidth]{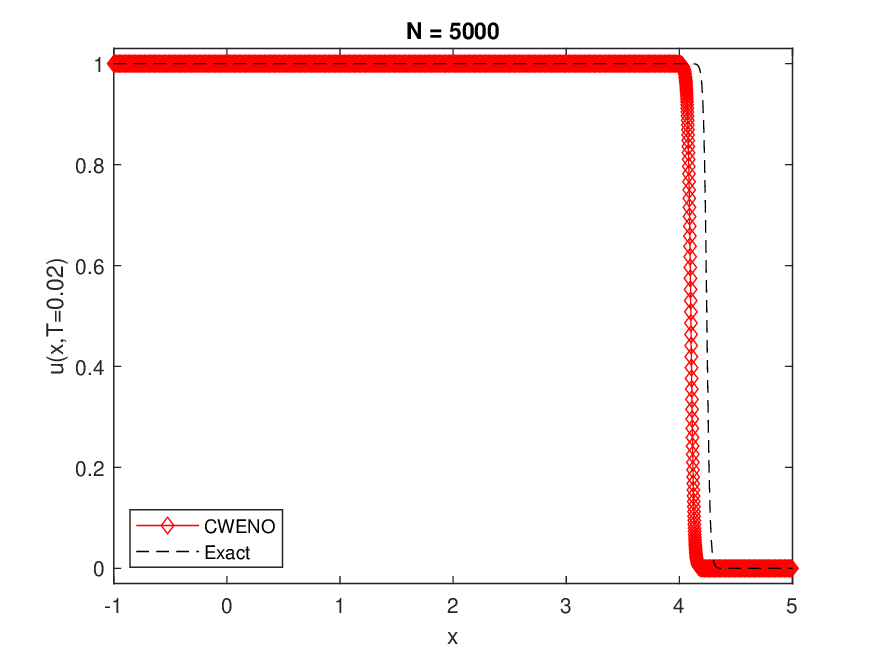}
\includegraphics[width=0.325\textwidth]{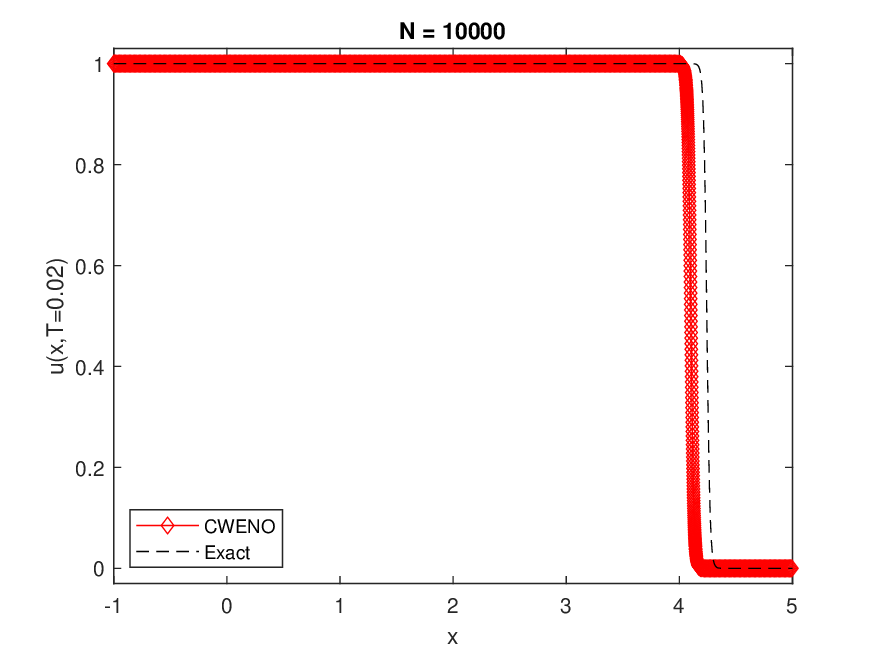}
\caption{Solution profiles for Newell-Whitehead-Segel equation ($\alpha=2$) at $T=0.02$ computed by CWENO with $N=1000$ (left), $5000$ (middle) and $10000$ (right).
The dashed black line is the exact solution.}
\label{fig:nws_alpha2_N}
\end{figure}

\begin{figure}[htbp]
\centering
\includegraphics[width=0.325\textwidth]{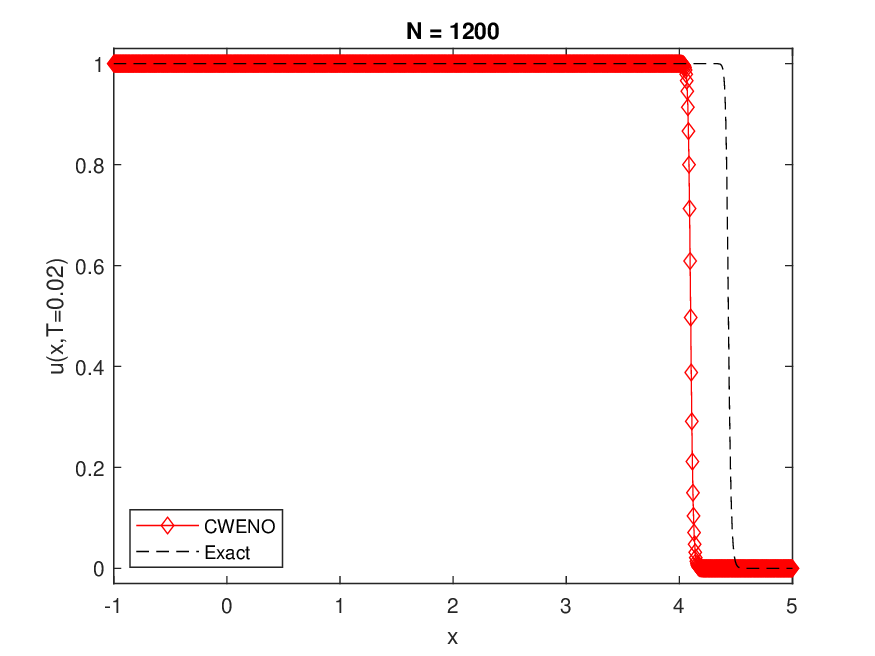}
\includegraphics[width=0.325\textwidth]{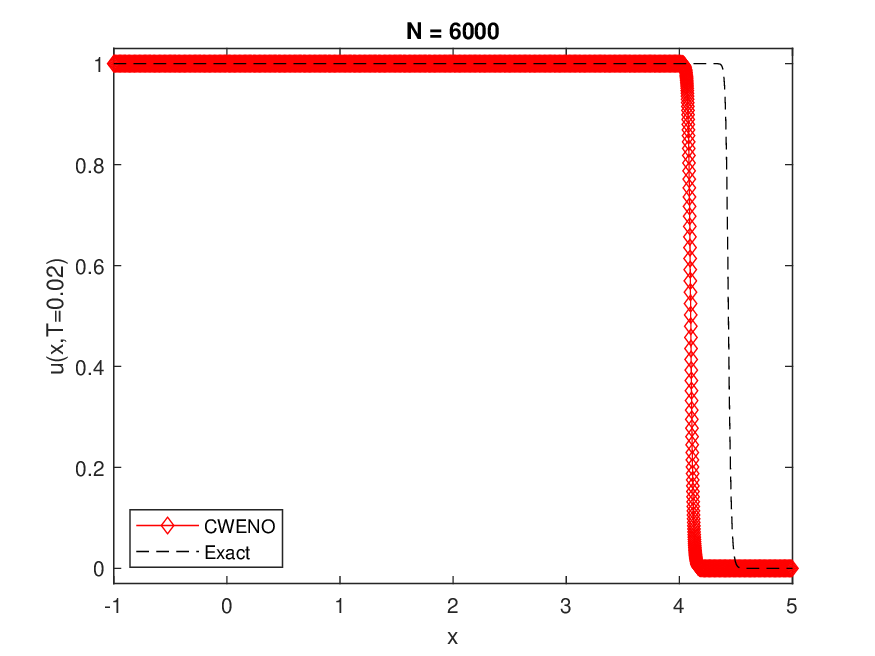}
\includegraphics[width=0.325\textwidth]{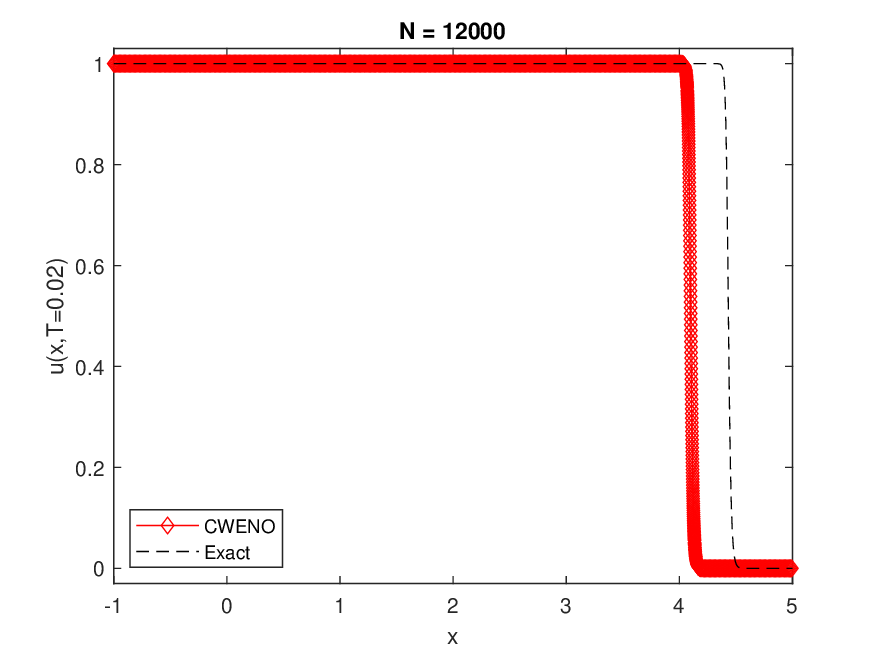}
\caption{Solution profiles for Newell-Whitehead-Segel equation ($\alpha=3$) at $T=0.02$ computed by CWENO with $N=1200$ (left), $6000$ (middle) and $12000$ (right).
The dashed black line is the exact solution.}
\label{fig:nws_alpha3_N}
\end{figure}

\begin{figure}[htbp]
\centering
\includegraphics[width=0.325\textwidth]{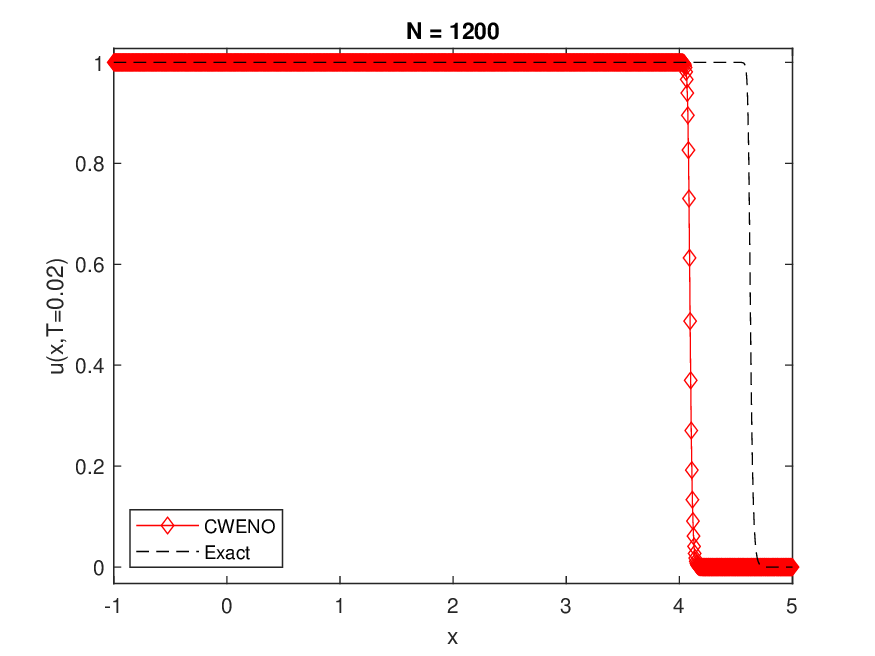}
\includegraphics[width=0.325\textwidth]{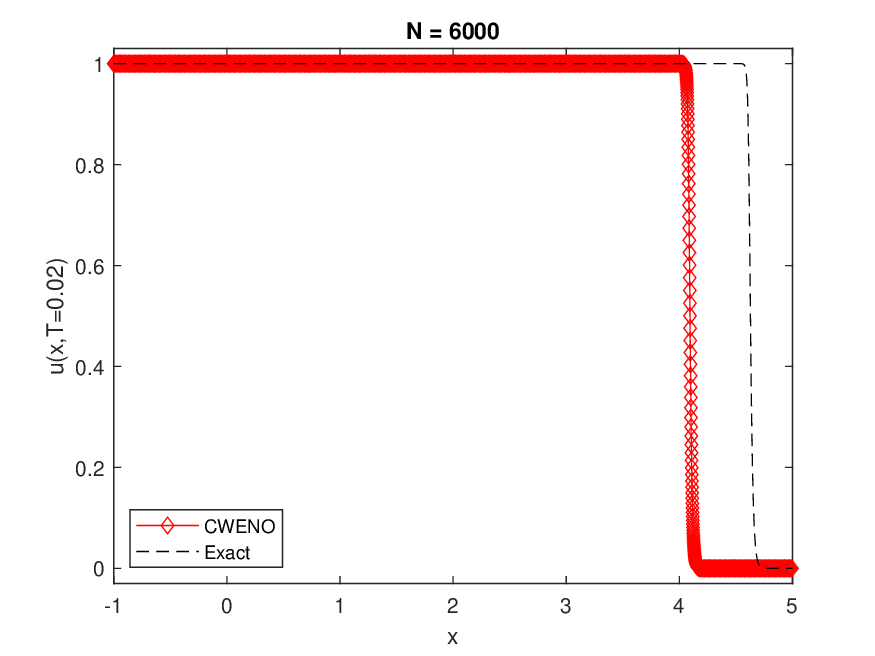}
\includegraphics[width=0.325\textwidth]{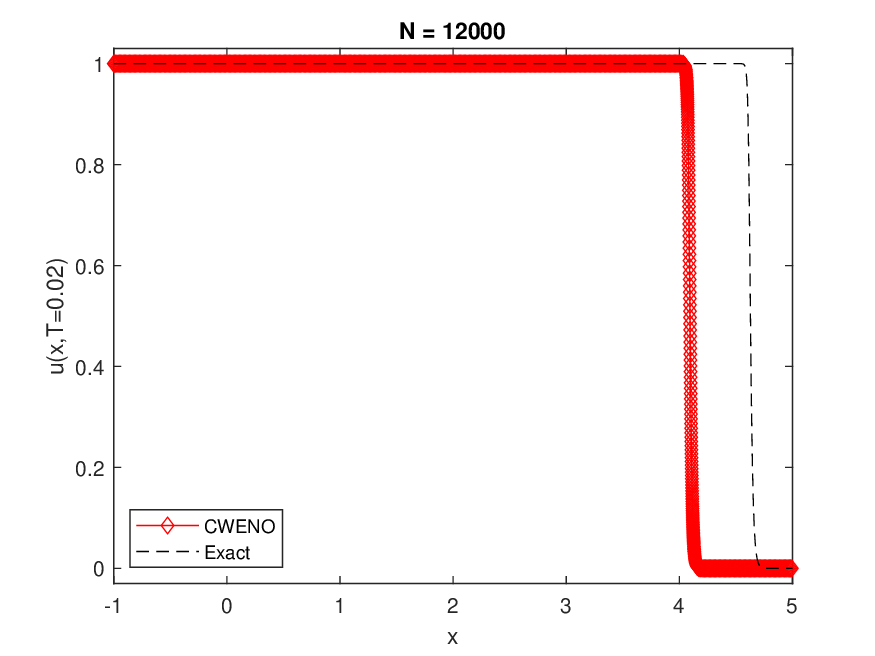}
\caption{Solution profiles for Newell-Whitehead-Segel equation ($\alpha=4$) at $T=0.02$ computed by CWENO with $N=1200$ (left), $6000$ (middle) and $12000$ (right).
The dashed black line is the exact solution.}
\label{fig:nws_alpha4_N}
\end{figure}

For $\alpha=2, \, 3$ and $4$, there is no dramatic improvement in reducing the lag in Figures \ref{fig:nws_alpha2_N}, \ref{fig:nws_alpha3_N} and \ref{fig:nws_alpha4_N}, respectively.
In Table \ref{tab:nws_error_N}, it shows that the errors are even increasing for $\alpha=4$ as we use the finer grid. 

\begin{table}[htbp]
\renewcommand{\arraystretch}{1.1}
\scriptsize
\centering
\caption{$L_1, \, L_2, \, L_\infty$ errors of CWENO for Newell-Whitehead-Segel equation.}      
\begin{tabular}{ccccc} 
\hline
$\alpha$ & $N$ & $L_1$ error & $L_2$ error & $L_\infty$ error \\ 
\hline 
\multirow{3}{*}{$1$} & $500$  & 0.010764 & 0.070489 & 0.662911 \\  
                     & $1000$ & 0.001107 & 0.007766 & 0.080663 \\  
                     & $2000$ & 0.000810 & 0.005687 & 0.059123 \\
\hline
\multirow{3}{*}{$2$} & $1000$  & 0.023976 & 0.139327 & 0.990892 \\  
                     & $5000$  & 0.023551 & 0.137774 & 0.989952 \\
                     & $10000$ & 0.023520 & 0.137661 & 0.989877 \\
\hline
\multirow{3}{*}{$3$} & $1200$  & 0.055693 & 0.227519 & 1-4.362e-07 \\  
                     & $6000$  & 0.055693 & 0.227511 & 1-4.426e-07 \\  
                     & $12000$ & 0.055677 & 0.227474 & 1-4.454e-07 \\
\hline
\multirow{3}{*}{$4$} & $1200$  & 0.089099 & 0.292279 & 1-1.956e-12 \\  
                     & $6000$  & 0.089238 & 0.292511 & 1-1.922e-12 \\  
                     & $12000$ & 0.089268 & 0.292563 & 1-1.905e-12 \\
\hline
\end{tabular}
\label{tab:nws_error_N}
\end{table}

\begin{figure}[htbp]
\centering
\includegraphics[width=0.335\textwidth]{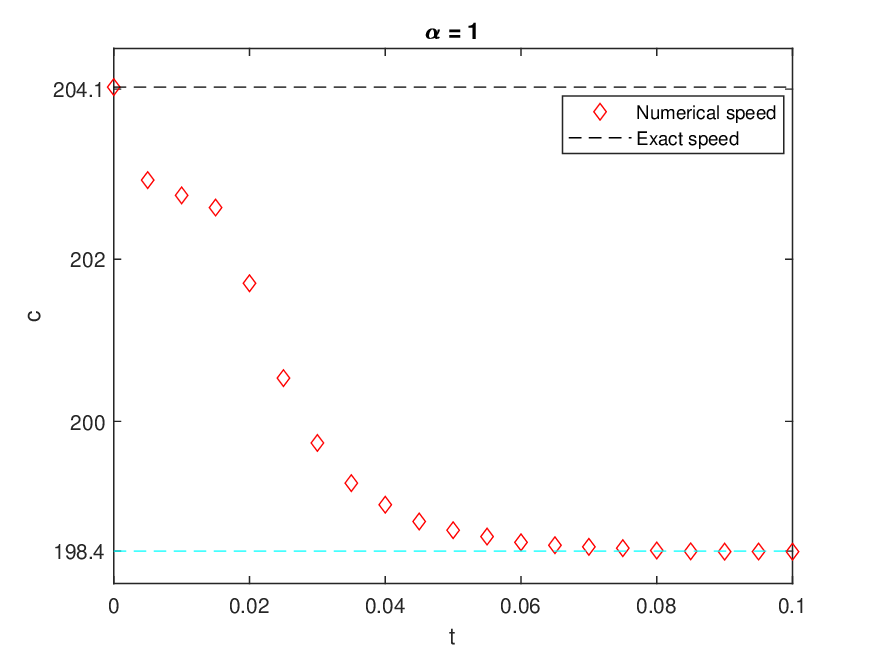}
\includegraphics[width=0.335\textwidth]{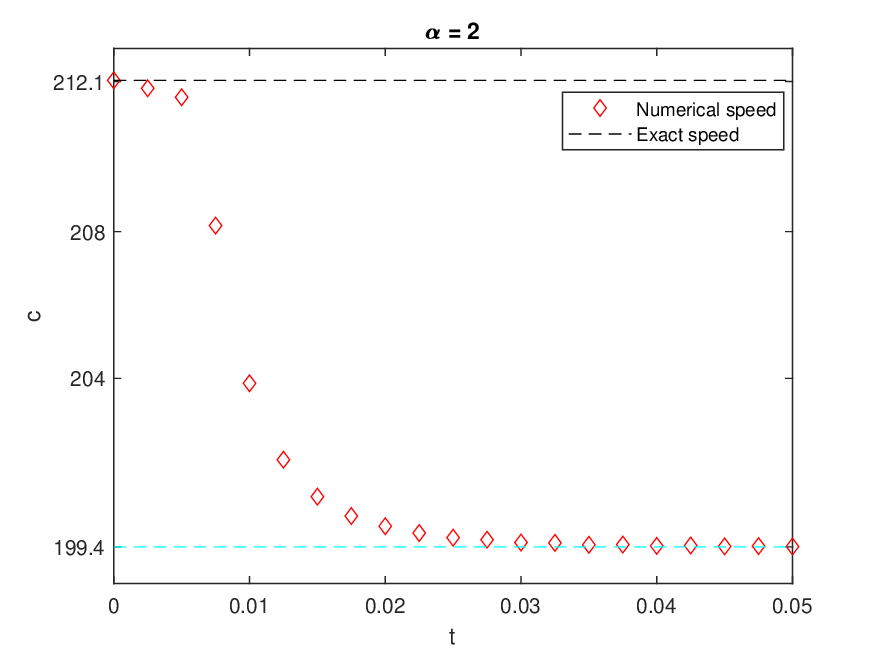}
\includegraphics[width=0.335\textwidth]{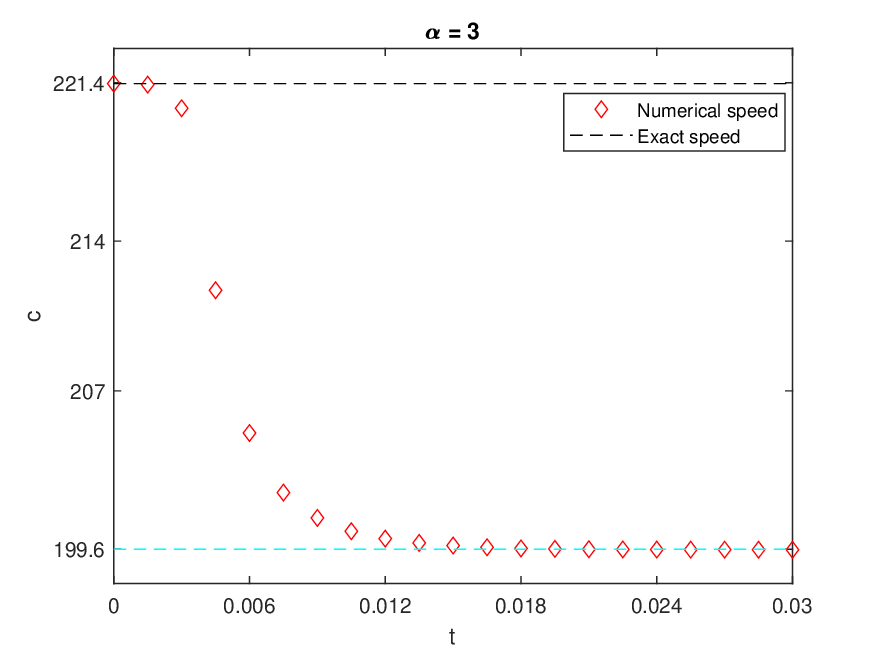}
\includegraphics[width=0.335\textwidth]{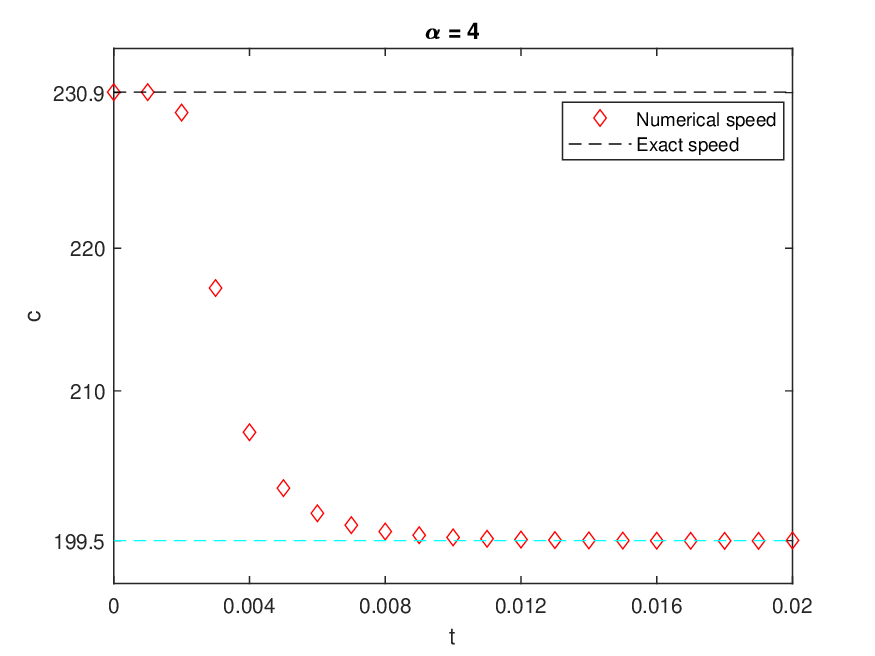}
\caption{Estimate of the numerical speed for Newell-Whitehead-Segel equation with $\alpha=1$ (top left), $2$ (top right), $3$ (bottom left) and $4$ (bottom right) by CWENO.
The dashed black line corresponds to the exact speed. The cyan line corresponds to the converged speed.}
\label{fig:nws_speed}
\end{figure}

In \cite{Gazdag}, it was observed that the numerical simulation of Fisher's equation is unstable, i.e., if the inital wave takes the velocity that is greater than the minimum speed, the initial wave could evolve into the wave of the minimum speed, due to the general perturbations of infinite extent. 
Figure \ref{fig:nws_speed} shows the estimate of the numerical speed of the sharp front for the Newell-Whitehead-Segel equation by CWENO for $\alpha=1, \, 2, \, 3, \, 4$. 
We can see that the numerical speed drops down and converges to a speed very close to the theoretical minimum given by $c_{\min} = 2 \sqrt{\rho D} = 200$, as obtained in \cite{MurrayI} for Fisher's equation.

\begin{figure}[htbp]
\centering
\includegraphics[width=0.325\textwidth]{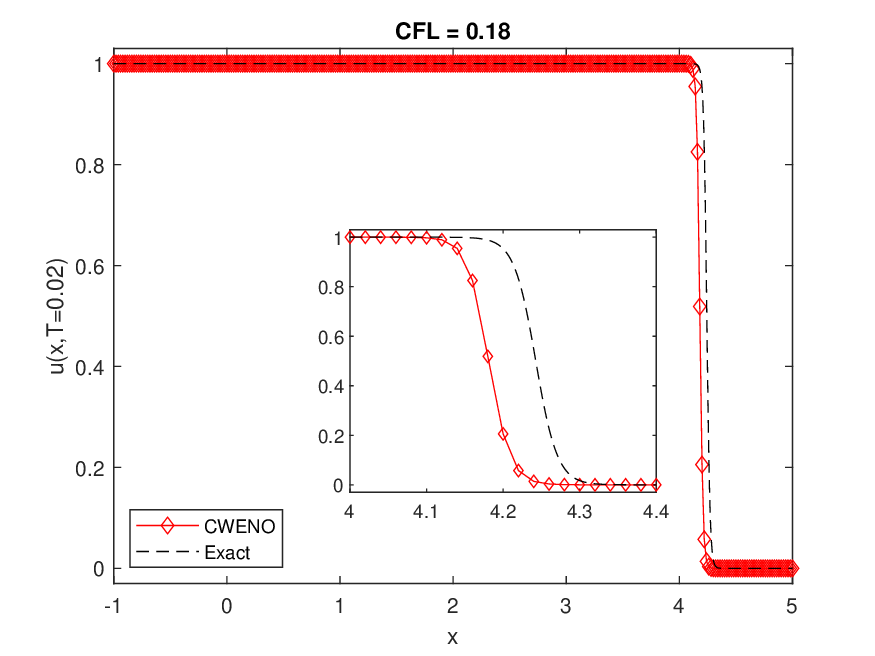}
\includegraphics[width=0.325\textwidth]{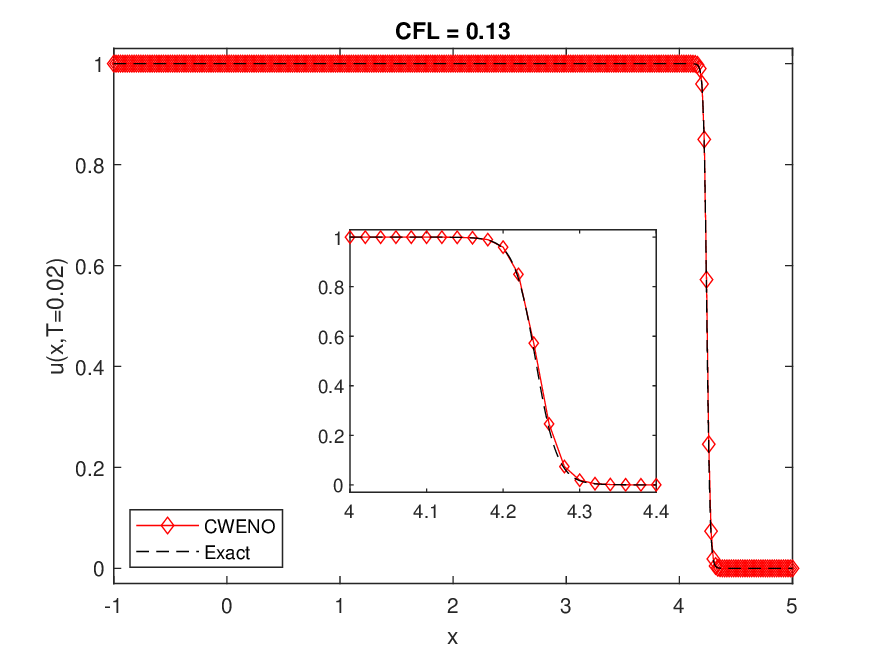}
\includegraphics[width=0.325\textwidth]{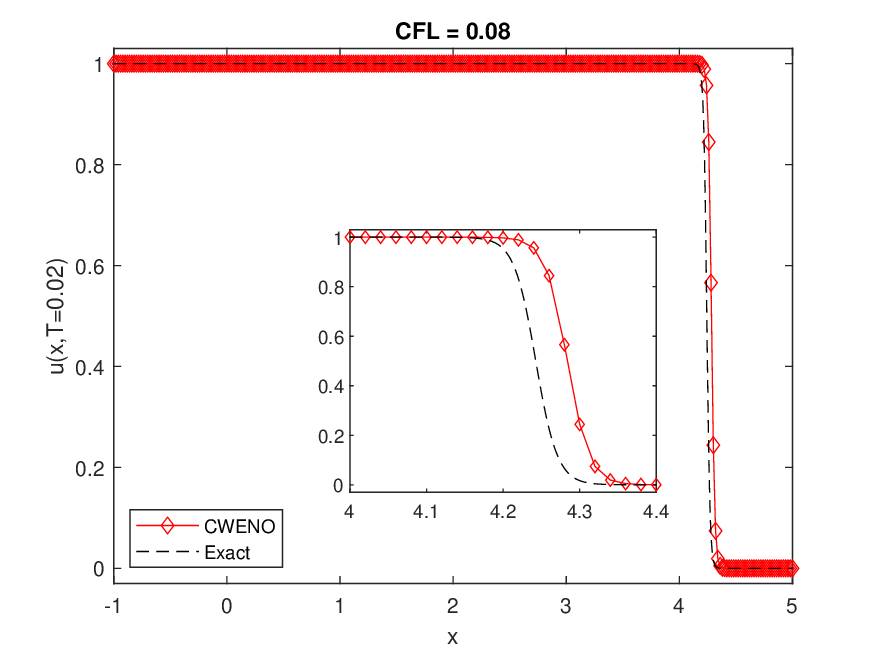}
\caption{Solution profiles for Newell-Whitehead-Segel equation ($\alpha=2$) at $T=0.02$ computed by central WENO (right) with $\cfl=0.18$ (left), $0.13$ (middle) and $0.08$ (right).
The dashed black line is the exact solution.}
\label{fig:nws_alpha2_CFL}
\end{figure}

\begin{figure}[htbp]
\centering
\includegraphics[width=0.325\textwidth]{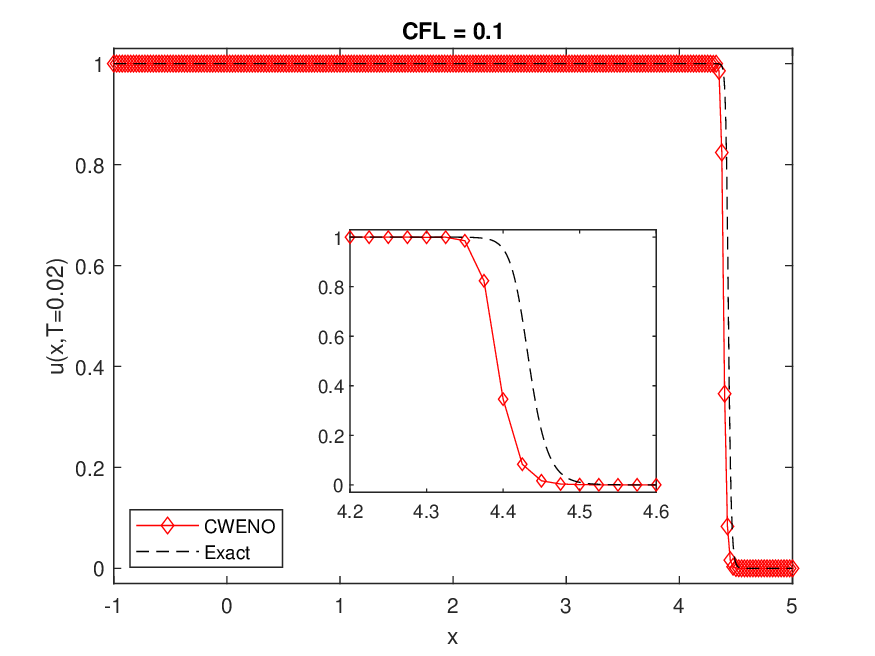}
\includegraphics[width=0.325\textwidth]{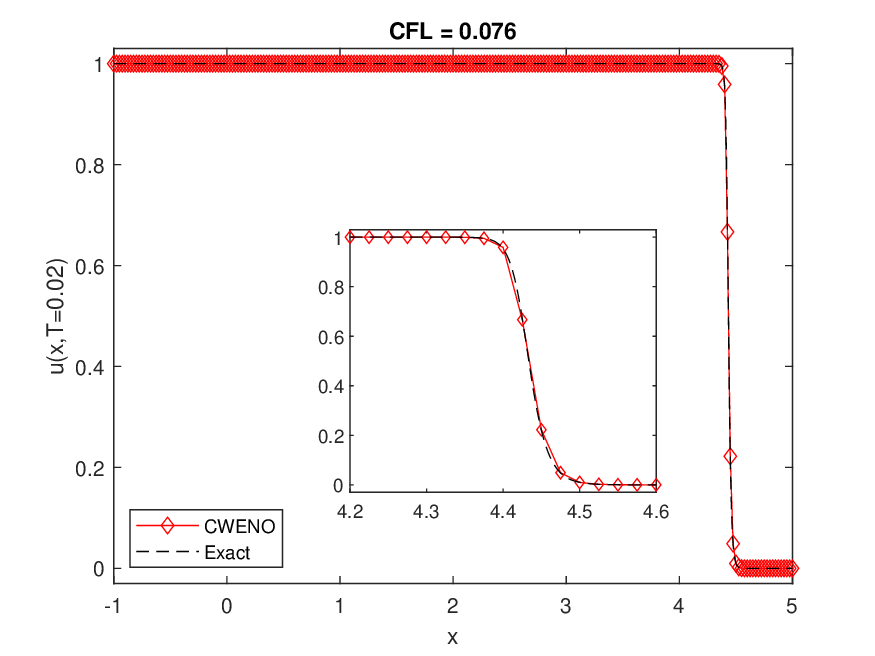}
\includegraphics[width=0.325\textwidth]{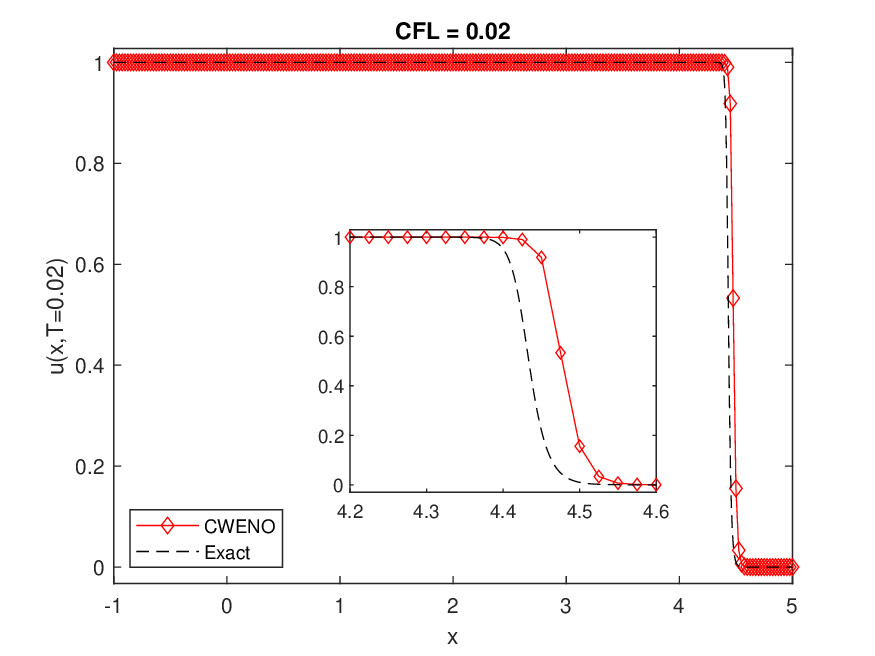}
\caption{Solution profiles for Newell-Whitehead-Segel equation ($\alpha=3$) at $T=0.02$ computed by central WENO (right) with $\cfl=0.1$ (left), $0.08$ (middle) and $0.02$ (right).
The dashed black line is the exact solution.}
\label{fig:nws_alpha3_CFL}
\end{figure}

\begin{figure}[htbp]
\centering
\includegraphics[width=0.325\textwidth]{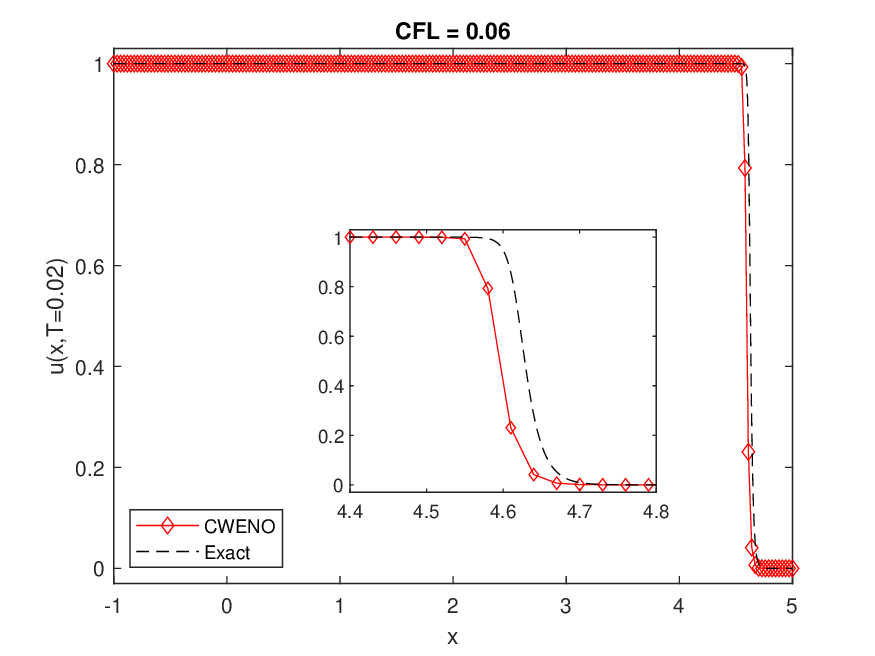}
\includegraphics[width=0.325\textwidth]{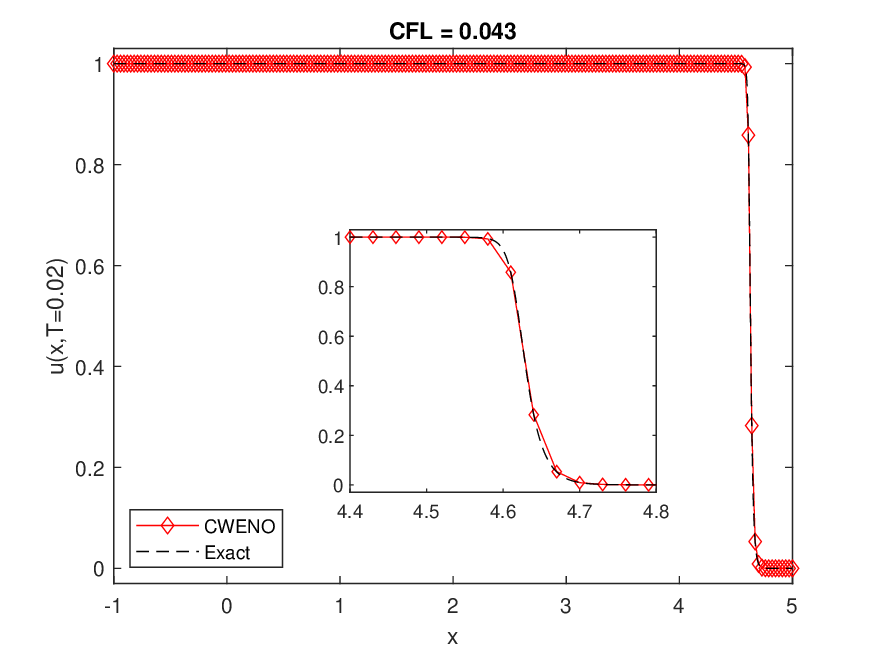}
\includegraphics[width=0.325\textwidth]{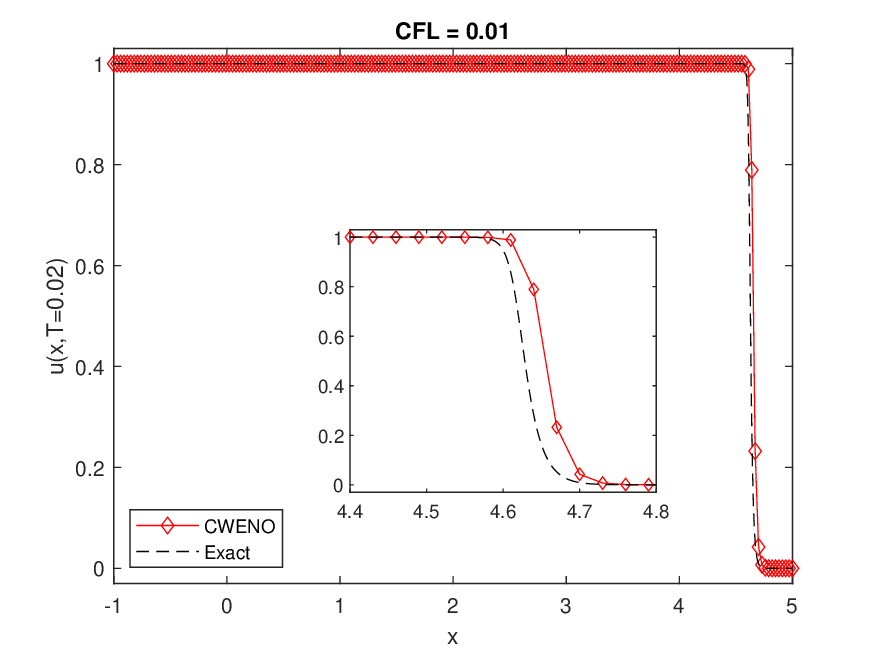}
\caption{Solution profiles for Newell-Whitehead-Segel equation ($\alpha=4$) at $T=0.02$ computed by central WENO (right) with $\cfl=0.06$ (left), $0.04$ (middle) and $0.01$ (right).
The dashed black line is the exact solution.}
\label{fig:nws_alpha4_CFL}
\end{figure}

From a large number of numerical experiments, we find that reducing both $N$ and the CFL number is one possible way to improve the accuracy of CWENO.
However, it is not an easy task to find such an optimal combination of $N$ and the CFL number, where the numerical speed is close to the exact one.
Figures \ref{fig:nws_alpha2_CFL}, \ref{fig:nws_alpha3_CFL} and \ref{fig:nws_alpha4_CFL} show that how decreasing both $N$ and CFL value makes the sharp front from delay to acceleration for the Newell–Whitehead–Segel equation.
The corresponding errors are displayed in Table \ref{tab:nws_error_N_CFL}, which shows how the numerical speed behaves compared to the exact speed for different values of $\alpha$ and CFL numbers. 
The table summarizes the results on the wave front speed from Figures \ref{fig:nws_alpha2_CFL}, \ref{fig:nws_alpha3_CFL} and \ref{fig:nws_alpha4_CFL}. 
As shown in the table, for the given $N$, decreasing CFL number does not necessarily yield convergence numerical speed. 
Instead, it seems that there exists optimal CFL number that provides the numerical speed which is similar to the exact speed. 
This phenomenon in this paper is not observed for linear problems and other equations considered rather than the Newell-Whitehead-Segel equation. 
At the moment, it is not known how to find the optimal CFL number for given $N$ and $\alpha$ for the Newell-Whitehead-Segel equation, which will be investigated in our future work. 

\begin{table}[htbp]
\renewcommand{\arraystretch}{1.1}
\scriptsize
\centering
\caption{$L_1, \, L_2$ and $L_\infty$ errors for Newell-Whitehead-Segel equation.}      
\begin{tabular}{ccccccc} 
\hline
$\alpha$ & $N$ & CFL & $L_1$ error & $L_2$ error & $L_\infty$ error & \\ 
\hline
$2$ & $300$ & $0.18$ & 0.010184 & 0.076233 & 0.774677 & $< c_\exact$ \\  
    &       & $0.13$ & 0.000270 & 0.002189 & 0.025498 & $\approx c_\exact$ \\
    &       & $0.08$ & 0.006861 & 0.054174 & 0.618005 & $> c_\exact$ \\
\hline
$3$ & $240$ & $0.1$   & 0.006816 & 0.057201 & 0.606482 & $< c_\exact$ \\  
    &       & $0.076$ & 0.000066 & 0.000590 & 0.006424 & $\approx c_\exact$ \\  
    &       & $0.02$  & 0.007212 & 0.059344 & 0.698676 & $> c_\exact$ \\
\hline
$4$ & $200$ & $0.06$  & 0.005626 & 0.049514 & 0.626126 & $< c_\exact$ \\  
    &       & $0.043$ & 0.000026 & 0.000164 & 0.001466 & $\approx c_\exact$ \\  
    &       & $0.01$  & 0.004301 & 0.039159 & 0.507090 & $> c_\exact$ \\
\hline
\end{tabular}
\label{tab:nws_error_N_CFL}
\end{table}

\section{Conclusions} \label{sec:conclusions}
In this work, we applied several sixth-order finite difference WENO methods to approximate travelling wave solutions in reaction-diffusion systems, including the Fisher's, Zeldovich, Newell-Whitehead-Segel and bistable equations, and the Lotka-Volterra competition-diffusion system. 
These equations are highly nonlinear and yield travelling wave solutions with sharp wave fronts. 
The popularly used numerical methods in the reaction-diffusion community is the high order finite difference methods. 
However, due to the sharp wave fronts, the finite difference methods require highly large number of grid point which strictly limits the time stepping. 
As the WENO method is a good alternative to the finite difference method for dealing with the sharp solution profile, we utilize the WENO methods developed for the parabolic term and demonstrate that the WENO, particularly, the central WENO method yields better performance than the commonly used finite difference method in the community and other WENO methods. 
That implies that with a smaller number of grid points and a large value of time stepping could be used for the same degree of accuracy by the finite difference method which should use much large number of grid points and smaller CFL number. 
Thus, it is highly recommended that the WENO method is used for the reaction-diffusion systems that are highly nonlinear and yield travelling wave solutions. 

For the travelling wave solution, approximating the wave speed properly is an important task. 
We find numerically that an accurate numerical wave speed is obtained by the optimal combinations of the number of grid points and the CFL number in the Newell-Whitehead-Segel equation. 
That is, for the given number of grid points, simply decreasing arbitrarily the CFL number does not yield an accurate wave speed. 
Instead, there exists an optimal value of CFL number that provides the numerical wave speed that converges to the exact speed. 
At the moment the way to find such an optimal CFL number is not known. 
In our future work we will investigate this phenomenon more rigorously and attempt to find the optimal value of CFL number with the given number of grid points.

\section*{Acknowledgments}
This work was supported by NRF of Korea under the grant number 2021R1A2C3009648 and POSTECH Basic Science Research Institute under the NRF grant number NRF2021R1A6A1A1004294412. 
JHJ is also supported partially by NRF grant by the Korea government (MSIT) (RS-2023-00219980).


\providecommand{\bysame}{\leavevmode\hbox to3em{\hrulefill}\thinspace}
\providecommand{\MR}{\relax\ifhmode\unskip\space\fi MR }
\providecommand{\MRhref}[2]{%
  \href{http://www.ams.org/mathscinet-getitem?mr=#1}{#2}
}
\providecommand{\href}[2]{#2}

\end{document}